\documentclass[11pt]{article}

\usepackage{amssymb,amsmath,latexsym, bbm}
\usepackage{epsfig}

\newtheorem{thm}{Theorem}[section]

\newtheorem{cor}[thm]{Corollary}

\newtheorem{lem}[thm]{Lemma}
\newtheorem{propn}[thm]{Proposition}

\newcommand{\pf}{\noindent{\bf Proof.} }
\def\qed{{\hfill $\Box$ \bigskip}}

\def\R{\mathbb{R}}
\def\bN{\mathbb{N}}
\def\bz{\mathbb{Z}}
\def\bP{\mathbb{P}}
\def\bQ{\mathbb{Q}}
\def\E{\mathbb{E}}

\def\sE{\mathcal{E}}
\def\sF{\mathcal{F}}

\def\NN{\mathcal{N}}

\def\LL{\mathcal{L}}

\def\lam{\lambda}

\def\1{\mathbbm{1}}

\def\bar{\overline}
\def\wh{\widehat}
\def\wt{\widetilde}
\def\eps{\varepsilon}
\def\tp {\wt \phi}
\def\angel#1{{\langle#1\rangle}}

\def\bee{\begin{equation}}
\def\eee{\end{equation}}

\makeatletter
\@addtoreset{equation}{section}

\makeatother

\begin{document}
\bibliographystyle{plain}

\title{\Large \bf
A priori H\"older estimate, parabolic Harnack principle and heat
kernel estimates   for diffusions with jumps}

\author{{\bf Zhen-Qing Chen}
\thanks{Research partially supported by NSF Grant  DMS-06000206.}
\quad and  \quad {\bf Takashi~Kumagai}
\thanks{Research partially supported by the
Grant-in-Aid for Scientific Research (B) 18340027.}
}
\date{(August 18, 2008)}
\maketitle

\begin{abstract} In this paper, we consider the following type of non-local
 (pseudo-differential) operators $\LL $ on $\R^d$:
 $$
  \LL  u(x) =\frac12 \sum_{i, j=1}^d \frac{\partial}{\partial x_i} \left(a_{ij}(x)
 \frac{\partial}{\partial x_j}\right) +
 \lim_{\eps \downarrow 0} \int_{\{y\in \R^d: \, |y-x|>\eps\}}
 (u(y)-u(x)) J(x, y) dy,
 $$
 where $A(x)=(a_{ij}(x))_{1\leq i,
j\leq d}$ is a measurable $d\times d$ matrix-valued function on $\R^d$ that is
uniform elliptic
and bounded and $J$ is a symmetric measurable
non-trivial non-negative kernel on $\R^d\times \R^d$ satisfying certain conditions.
Corresponding to $\LL$ is a symmetric strong Markov process $X$ on $\R^d$ that has
both the diffusion component and pure jump component. We establish a priori H\"older
estimate for bounded parabolic functions of $\LL$ and parabolic Harnack principle
for positive parabolic functions of $\LL$. Moreover, two-sided sharp heat kernel
estimates are derived for such operator $\LL$ and jump-diffusion $X$.
In particular, our results apply to the mixture of symmetric diffusion of uniformly
elliptic divergence form operator and   mixed stable-like processes on $\R^d$.
To establish these results, we employ methods from both probability theory and
analysis.
\end{abstract}

\section{Introduction}

 It is well-known that there is an intimate interplay between
 self-adjoint pseudo-differential  operators on $\R^d$
 and symmetric strong Markov processes on $\R^d$.
 For a large class of self-adjoint pseudo-differential  operators $\LL$ on $\R^d$
 that enjoys maximum property,
   there is a jump-diffusion $X$ on $\R^d$
associated with it so that $\LL$ is the infinitesimal generator of $X$, and vice
versa.  The connection between $\LL$ and $X$ can also be seen as
follows.
  The fundamental solution (also called heat kernel)
  for $\LL$ is the transition density function of $X$.
In this paper,  we are interested in the a priori H\"older estimate
for harmonic functions of such operator $\LL$, parabolic Harnack principle
and the sharp estimates on the heat
kernel of $\LL$.

Throughout this paper,
$d\geq 1$ is an integer.
Denote by $m_d$ the $d$-dimensional Lebesgue measure in $\R^d$, and
$C^1_c(\R^d)$ the space of $C^1$-functions on $\R^d$ with compact
support.
  We consider the following type of non-local
 (pseudo-differential) operators $\LL $ on $\R^d$:
 \begin{equation}\label{e:op}
  \LL  u(x) =\frac12 \sum_{i, j=1}^d \frac{\partial}{\partial x_i} \left(a_{ij}(x)
 \frac{\partial}{\partial x_j}\right) +
 \lim_{\eps \downarrow 0} \int_{\{y\in \R^d: \, |y-x|>\eps\}}
 (u(y)-u(x)) J(x, y) dy,
 \end{equation}
where $A(x)=(a_{ij}(x))_{1\leq i,
j\leq d}$ is a measurable $d\times d$ matrix-valued function on $\R^d$ that is
uniform elliptic and
bounded in the sense that there exists a constant $c\geq 1$ such that
\begin{equation}\label{unielli}
c^{-1}\sum_{i=1}^d\xi_i^2 \le \sum_{i,j=1}^d a_{ij}(x)\xi_i\xi_j\le
c\sum_{i=1}^d\xi_i^2 \qquad \hbox{for every  } x,
(\xi_1,\cdots,\xi_d)\in \R^d,
\end{equation}
and $J$ is a symmetric non-negative measurable kernel on $\R^d\times \R^d$
  such that there are positive constants $\kappa_0>0 $,
 and $\beta \in (0, 2)$ so that
\begin{equation}\label{e:J1}
  J(x, y) \leq \kappa_0
|x-y|^{-d-\beta } \qquad \hbox{for } |x-y|
 \leq \delta_0,
 \end{equation}
and that
\begin{equation}\label{e:J2}
\sup_{x\in \R^d} \int_{\R^d} (|x-y|^2 \wedge 1) J(x, y)\, dy <\infty .
\end{equation}
 Clearly under condition \eqref{e:J1}, condition \eqref{e:J2} is equivalent to
$$ \sup_{x\in \R^d} \int_{\{y\in \R^d: |y-x|\geq 1\}}  J(x, y) \, dy <\infty .
$$

 Associated with such a non-local operator $\LL$
 is an $\R^d$-valued symmetric strong Markov process $X$ whose associated
 Dirichlet form  $(\sE, \sF)$
  on $L^2(\R^d; m_d)$ is given by
 \begin{equation}\label{e:DF}
\begin{cases} \displaystyle
\sE (u,v)=\frac12 \int_{\R^d} \nabla u(x)\cdot A(x) \nabla
v(x)dx+\int_{\mathbb R^d}
(u(x)-u(y))(v(x)-v(y))J(x,y)dxdy , \\
\hskip 0.3truein \sF \, = \, \overline{C^1_c (\R^d)}^{\sE_1} ,
\end{cases}
\end{equation}
where for $\alpha >0$, $\sE_\alpha (u, v):= \sE(u, v) + \alpha
\int_{\R^d} u(x) v(x) m_d (dx)$.

When the jumping kernel $J\equiv 0$ in \eqref{e:op} and \eqref{e:DF},
$\LL$ is a uniform elliptic operator of divergence form and $X$
is a symmetric diffusion on $\R^d$. It is well-known that
$X$ has a joint H\"older continuous transition density function $p(t, x, y)$,
which enjoys the following celebrated
Aronson's two-sided heat kernel estimate:
there are constants $c_k>0$,
 $k=1, \cdots, 4$, so that
$$ c_1 \, p^c(t, c_2|x-y|) \leq p(t, x, y)\leq c_3\, p^c(t, c_4|x-y|) \qquad \hbox{for }
t>0, x, y \in \R^d.
 $$
 Here
 \begin{equation}\label{e:pc}
  p^c(t, r) := t^{-d/2} \exp ( - r^2/t).
 \end{equation}
It is also known that parabolic Harnack principle holds for such $\LL$ and that
every  bounded parabolic function of $\LL$ is locally H\"older continuous.
See \cite{Str} for some history and a survey on this subject, where a mixture of analytic and probabilistic method is presented.

 Let $\phi$ be a strictly increasing continuous
function $\phi: \R_+\to \R_+$ with $\phi (0)=0$, and $\phi(1)=1$
such that there are constants $c\ge 1$, $0<\beta_1\leq \beta_2 <2$
such that
\begin{equation}\label{polycon}
c^{-1} \Big(\frac Rr\Big)^{\beta_1} \leq   \frac{\phi (R)}{\phi (r)}
\leq c \, \Big(\frac Rr\Big)^{\beta_2} \qquad \hbox{for every }
0<r<R<\infty,
\end{equation}
and
\begin{equation}\label{polycon2}
\int_0^r \frac{s}{\phi (s)} ds \leq c \, \frac{r^2}{\phi (r)} \qquad
\hbox{for every } r>0.
\end{equation}
Observe that  condition \eqref{polycon}  implies that
$$ c^{-1} r^{\beta_1} \leq \phi (r) \leq c r^{\beta_2}
\qquad \hbox{for } r\geq 1
$$
and
$$ c^{-1} r^{\beta_2} \leq \phi (r) \leq c r^{\beta_1}
\qquad \hbox{for } r\in (0, 1].
$$
In the sequel,
if $f$ and $g$ are two functions defined on a set $D$, $f\asymp g$
means that there exists $C>0$ such that $C^{-1}f(x)\le g(x)\le
C\,f(x)$ for all $x\in D$.

 When $A(x)\equiv 0$ in \eqref{e:DF} and $J$ is given by
  \begin{equation}\label{e:J4}
J(x, y) \asymp \frac1{|x-y|^{d}\, \phi (|x-y|)},
\end{equation}
where $\phi$ satisfies the conditions \eqref{polycon}-\eqref{polycon2},
the corresponding process $X$ is a mixed stable-like process on $\R^d$
studied in \cite{CK2}.
A typical example of $J$ satisfying condition \eqref{e:J4} is
$$ J(x, y)= \int_{\alpha_1}^{\alpha_2} \frac{c(\alpha, x, y)}{|x-y|^{d+\alpha}}
 \, \nu (d\alpha),
 $$
 where $\nu$ is a probability measure on $[\alpha_1, \alpha_2]\subset (0, 2)$
 and $c(\alpha, x, y)$ is a symmetric function in $x$ and $y$ is bounded
 between two positive constants that are independent of $\alpha\in [\alpha_1, \alpha_2]$.
Under the above condition, a priori H\"older estimate and parabolic Harnack principle
are established in \cite{CK2} for parabolic functions of $X$. Moreover, it is proved
in \cite{CK2} that $X$ has a jointly continuous transition density function
$p(t, x, y)$ and that it has the following two-sided sharp estimates: there are
positive constants $0<c_1<c_2$ so  that
$$ c_1 p^j(t,  |x-y|) \leq p(t, x, y)\leq c_2 p^j(t,  |x-y|) \qquad \hbox{for }
t>0, x, y \in \R^d,
 $$
where
\begin{equation}\label{eqn:4}
 p^j(t, r) :=\left( \phi^{-1}(t)^{-d} \wedge
\frac{t}{r^d \phi (r)} \right)
\end{equation}
with $\phi^{-1}$ being the inverse function of $\phi$.
Here and in the sequel,   for two real numbers $a$ and $b$, $a\wedge
b:= \min \{a, b\}$ and $a\vee b:=\max \{a, b\}$.
We point out that, in contrast to the diffusions (or differential operator) case,
heat kernel estimates for pure jump processes (or non-local integro-differential
operators) have been studied only quite recently. See
the introduction part of
\cite{CK2} for a brief account
of some history.

In this paper, we consider the case where both $A$ and $J$ are non-trivial in
\eqref{e:op} and \eqref{e:DF}.  Clearly the corresponding operators and
jump diffusions take up an important place both in
 theory and in applications. However there are very limited
 work in literature for this mixture case
 on the topics of this paper,
 see \cite{BKU}, \cite{CKS} and \cite{SV07} though.
 One of the difficulties in obtaining fine properties for such an operator
 $\LL$ and process $X$ is that it exhibits different scales: the diffusion part
 has Brownian scaling $r\mapsto   r^2$ while the pure jump part has
 a different type of scaling.
 Nevertheless, there is a folklore which says that with the presence of
the diffusion part corresponding to $\frac12
\sum_{i, j=1}^d \frac{\partial}{\partial x_i} \left(a_{ij}(x)
 \frac{\partial}{\partial x_j}\right)$, better results can be expected
 under weaker assumptions on the jumping kernel $J$ as the diffusion
 part helps to smooth things out. Our investigation confirms such an
 intuition. In fact we can establish a priori H\"older estimate and
 parabolic Harnack inequality under weaker conditions than \eqref{e:J4}.
 We now present the main results of this paper.
Let $W^{1,2}(\R^d)$ denote the Sobolev space of order $(1, 2)$ on $\R^d$;
that is, $W^{1,2}(\R^d):=\{ f\in L^2(\R^d; m_d): \, \nabla f \in
 L^2(\R^d; m_d)\}$.
It is not difficult to show the following.

\begin{propn}\label{P:charaF}
Under the conditions
\eqref{unielli}-\eqref{e:J2},
the domain of the Dirichlet form of \eqref{e:DF} is characterized by
$$\sF= W^{1,2}(\R^d) =\{f\in L^2(\R^d; m_d): \sE(f,f)<\infty\}.
$$
\end{propn}

Let $X$ be the symmetric Hunt process on $\R^d$ associated with the
regular Dirichlet form $(\sE, \sF)$. It will be shown in Theorem \ref{consernew}
below that $X$ has infinite lifetime. Let
  $Z=\{Z_t:=(V_0-t, X_t), t\geq 0\}$
 denote the
space-time process of $X$.
We say that  a non-negative real valued Borel  measurable function
$h(t,x)$ on $[0, \infty)\times \R^d$ is {\it parabolic}
(or {\it caloric})
on $D=(a,b)\times B(x_0,r)$
if there is a properly exceptional set $\NN\subset \R^d$ such that
for every relatively compact open subset $D_1$ of $D$,
$$
h(t, x)=\E^{(t,x)} [h (Z_{\tau_{D_1}})]
$$
for every $(t, x)\in D_1\cap ([0,\infty)\times
  (\R^d\setminus \NN))$,
where
$\tau_{D_1}=\inf\{s> 0: \, Z_s\notin D_1\}$.
We remark that in \cite{CK1, CK2}   the space-time process is defined
to be $(V_0+t, X_t)$ but this is merely a notational difference.
In this paper, we first show
that any parabolic function of $X$ is H\"older continuous.
Recall that $\delta_0$ is the positive constant in condition
\eqref{e:J1}.

\begin{thm}\label{T:holder}
Assume that
the Dirichlet form $(\sE, \sF)$ given by \eqref{e:DF}
satisfies the conditions \eqref{unielli}-\eqref{e:J2}
and that
for every $0<r<\delta_0$,
\begin{equation}\label{e:J1b}
\inf_{x_0, y_0\in \R^d \atop |x_0-y_0|=r} \,
 \inf_{x\in B(x_0, \, r/16)} \int_{B(y_0, \, r/16)} J(x,z )dz>0.
 \end{equation}
Then for every $R_0\in (0,1]$,
there are constants $c=c(R_0)>0$ and $\kappa>0$ such that for
every $0<R\le R_0$ and every bounded parabolic function
$h$ in $Q(0, x_0, 2R):=(0, 4R^2) \times B(x_0, 2R)$,
\begin{equation}\label{holder1}
|h(s, x) -h(t, y)| \leq c \,  \| h \|_{\infty, R} \,
R^{-\kappa} \, \left(  |t-s|^{1/2} + |x-y|
\right)^\kappa
\end{equation}
holds for $(s, x), \, (t, y)\in Q(0, x_0, R)$,
where $\|h\|_{\infty, R}:=\sup_{(t,y)\in [0, \, 4R^2 ]
\times \R^d \setminus \NN } |h(t,y)|$.
In particular, $X$ has a jointly continuous  transition density function $p (t, x, y)$
with respect to the Lebesgue measure. Moreover, for every $t_0\in (0, 1)$
there are constants $c >0$ and $\kappa>0$ such that
  for any $t, \, s \in (t_0, \,  1]$
and $(x_i, y_i)\in
\R^d\times \R^d$ with $i=1, 2$,
\begin{equation}\label{eqn:holder2}
 |p (s, x_1, y_1) -p (t, x_2, y_2)| \leq c \,
  t_0^{-(d+\kappa)/2}
\left(  |t-s|^{1/2} + |x_1-x_2|+|y_1-y_2| \right)^\kappa.
\end{equation}
\end{thm}

In addition to \eqref{unielli}-\eqref{e:J2} and \eqref{e:J1b}, if there is a constant $c>0$
 such that
\begin{equation}\label{e:J3}
 J(x,y) \le \frac{c}{r^d}\int_{B(x,r)}
  J(z, y) dz
 \quad
\hbox { whenever $r\le \frac 12 |x-y|\wedge 1$},\, x,y\in \R^d,
\end{equation}
we show that the parabolic Harnack principle holds for non-negative
parabolic functions of $X$.
(Note that (\ref{e:J3}) was introduced in
\cite{BBK, CKK} and it was denoted as (UJS)${}_{\leq 1}$ there.)

\begin{thm}\label{T:PHI}
Suppose that
the Dirichlet form $(\sE, \sF)$ given by \eqref{e:DF}
satisfies the condition \eqref{unielli}-\eqref{e:J2}, \eqref{e:J1b} and \eqref{e:J3}.
For every $\delta\in (0, 1)$, there
exist constants $c_1=c_1( \delta)$ and $ c_2=c_2( \delta)>0$
such that for every $z\in \R^d$, $t_0\ge 0$,
$0<R\leq c_1$
and every non-negative function $u$ on $[0, \infty)\times \R^d$
that is parabolic
on $(t_0,t_0+6\delta R^2 )\times B(z,4R)$,
\begin{equation}\label{eqn:4.1}
\sup_{(t_1,y_1)\in Q_-}u(t_1,y_1)
\le c_2 \,
\inf_{(t_2,y_2)\in Q_+}u(t_2,y_2),
\end{equation}
where $Q_-=(t_0+\delta R^2,t_0+2\delta R^2)\times B(x_0,R)$ and
$Q_+=(t_0+3\delta R^2,t_0+ 4\delta R^2)\times B(x_0,R)$.
\end{thm}

Note that
elliptic versions of Theorem \ref{T:holder} and \ref{T:PHI}
are claimed in \cite{Foo} under similar assumptions,
however we have some difficulty to follow some of the arguments there.
 Clearly, our theorems
imply the elliptic versions given in \cite{Foo}.

We next derive two-sided heat kernel estimate for $X$ when $J(x, y)$
satisfies the condition \eqref{e:J4}. Clearly
\eqref{e:J1}-\eqref{e:J2}, \eqref{e:J1b}
 and \eqref{e:J3} are satisfied when
\eqref{e:J4} holds. Recall that functions $p^c(t, x, y)$ and $p^j(t, x, y)$
are defined by \eqref{e:pc} and \eqref{eqn:4}, respectively.

\begin{thm}\label{mainHK}
Suppose that \eqref{unielli} holds and that the jumping kernel $J$ of
the Dirichlet form $(\sE, \sF)$ given by \eqref{e:DF}
satisfies the condition \eqref{e:J4}.
Denote by $p(t, x, y)$ the continuous transition density function of  the
symmetric Hunt process $X$ associated with the regular Dirichlet
form $(\sE, \sF)$ of \eqref{e:DF} with the jumping kernel $J$ given
by \eqref{e:J4}. There are positive constants $c_i$, $i=1, 2, 3, 4$
such that for every $t>0$ and $x, y \in \R^d$,
\begin{eqnarray}
&& c_1\,    \left( t^{-d/2}\wedge\phi^{-1}(t)^{-d} \right)   \wedge
\left( p^c(t, c_2|x- y|)+p^j(t, |x- y|) \right)
\nonumber\\
&\leq & p(t, x, y) \leq c_3\, \left( t^{-d/2}\wedge\phi^{-1}(t)^{-d}
\right) \wedge \left( p^c(t, c_4 |x- y|)+p^j(t, |x- y|) \right).
\label{eq:HKj+dul}\end{eqnarray}
\end{thm}

The following figure shows which term is the dominant term in each region when
$\phi$ in \eqref{e:J4} is given by
$\phi(r)=r^\alpha$ with $0<\alpha<2$. It is worth mentioning that there is a
short-time short-distance region in $t\le R^2\le 1$ where the jump part is the
dominant term.

\centerline{\epsfig{file=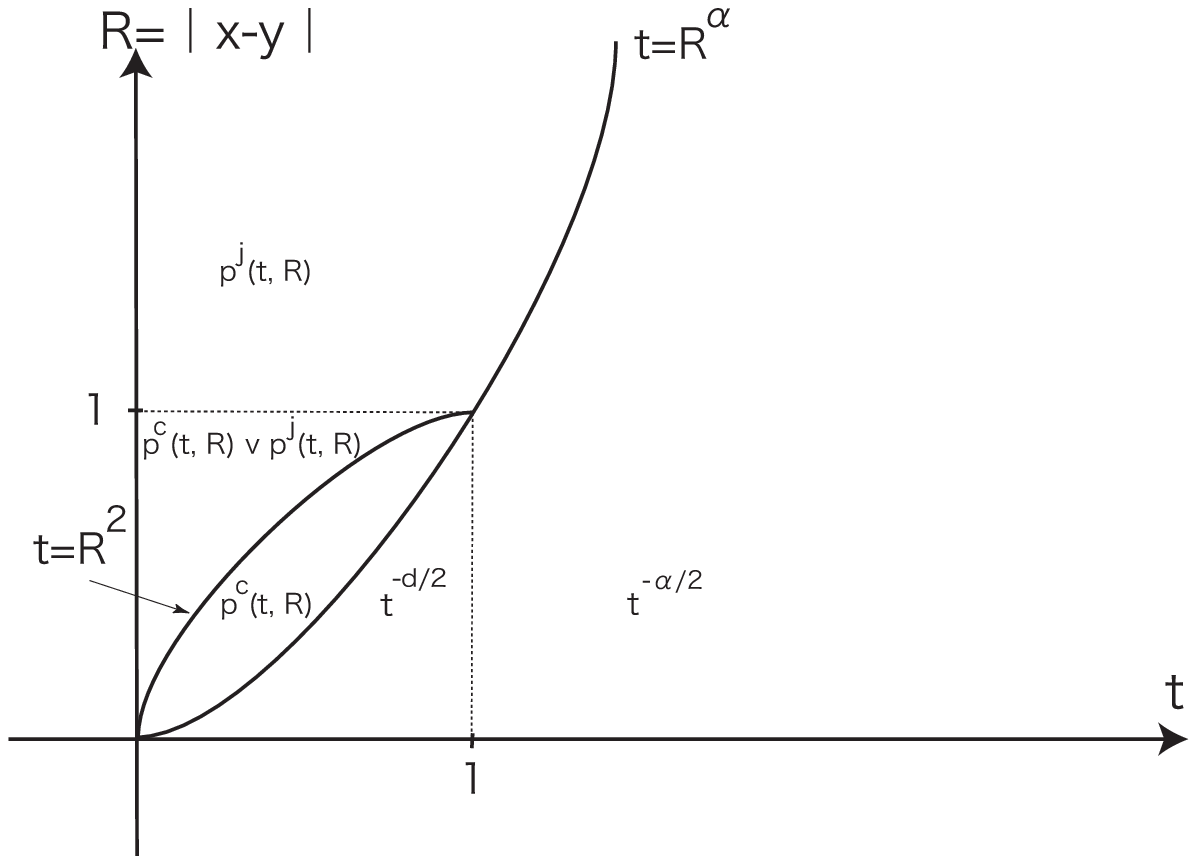, height=2in}}

When $A(x)\equiv I_{d\times d}$, the $d\times d$ identity matrix,  and
$J(x, y)= c |x-y|^{-d-\alpha}$ for some $\alpha \in (0, 2)$ in \eqref{e:DF},
that is, when $X$ is the independent sum of a Brownian motion $W$ on $\R^d$
and an isotropically symmetric $\alpha$-stable process $Y$ on $\R^d$,
the transition density function $p(t, x, y)$ can be expressed as the convolution
of the transition density functions of $W$ and $Y$,
 whose two-sided estimates are known.
In \cite{SV07},  heat kernel estimates for this L\'evy process $X$
are carried out by computing the convolution and the estimates  are given
in a form that depends on which region the point $(t, x, y)$ falls into.
Subsequently, the parabolic Harnack inequality \eqref{eqn:4.1}
 for such a L\'evy process
$X$ is derived in \cite{SV07} by using the two-sided Heat kernel estimate.
Clearly such an approach is not applicable in our setting even when $\phi (r)=r^\alpha$,
  since in our case, the diffusion and jumping part of $X$ are typically
not independent. The two-sided estimate in this simple form of
(\ref{eq:HKj+dul})
is a new observation even in the independent sum of a
Brownian motion and an isotropically symmetric $\alpha$-stable process case
considered in \cite{SV07}.

Our approach employs methods from both probability theory and analysis, but
it is mainly probabilistic. It   uses some ideas previously developed
in \cite{BBCK, BBK, CK1, CK2, CKK}.
To get a priori H\"older estimates for parabolic functions of $X$,
we establish the following three key ingredients.
\begin{description}
\item{(i)} Exit time upper bound estimate (Lemma \ref{L:exit}):
$$\E_x [ \tau_{B(x_0, r)}] \leq c_1r^2 \qquad \hbox{for } x\in B(x_0, r),
 $$
 where $\tau_{B(x_0, r)}:=\inf\{t>0: X_t\notin B(x_0, r)\}$
  is the first exit time from $B(x_0, r)$ by $X$.

\item{(ii)} Hitting probability estimate (\eqref{e:hitting} below):
$$
\bP_x \left(X_{\tau_{B(x, r)}} \notin B(x, s) \right)
\leq  \frac{c_2r^2}{(s\wedge 1)^2} \qquad \hbox{for every } r\in (0, 1]
\hbox{ and } s\geq 2r.
$$

\item{(iii)} Hitting probability estimate for space-time process $Z_t=(V_0-t, X_t)$
(Lemma \ref{L:exit3}):
 for every $x\in \R^d$, $r\in (0, 1]$ and any compact subset $A\subset Q(x, r):=(0, r^2)\times B(x, r)$,
 $$
  \bP^{(r^2, x)}
 (\sigma_A <\tau_r) \geq c_3
 \frac{m_{d+1}(A)}{r^{d+2}},
 $$
 where by slightly abusing the notation,
 $\sigma_A:=\{t>0: Z_t\in A\}$ is the first hitting time of $A$, $\tau_r:=\inf\{t>0: Z_t\notin Q(x, r)\}$
  is the first exit time from $Q(x, r)$ by $Z$
  and $m_{d+1}$ is the
 Lebesgue measure on $\R^{d+1}$.
\end{description}

Throughout this paper, we use the following notations.
The probability law of the process $X$ starting from $x$ is denoted as $\bP_x$
and the mathematical expectation under it is denoted as $\E_x$, while
probability law of the space-time process $Z=(V, X)$ starting from
$(t, x)$, i.e. $(V_0,X_0)=(t,x)$,
is denoted as $\bP^{(t, x)}$
and the mathematical expectation under it is denoted as $\E^{(t, x)}$.
To establish parabolic Harnack inequality, we need in addition the following.
\begin{description}\item{(iv)} Short time near-diagonal heat kernel estimate
(Theorem  \ref{T:3.4}): for every $t_0>0$, there is $c_4=c_4(t_0)>0$ such that
for every $x_0\in \R^d$ and $t\in (0, t_0]$,
$$ p^{B(x_0, \sqrt{t})} (t, x, y) \geq c_4 t^{-d/2} \qquad \hbox{for }
x, y\in B(x_0, \sqrt{t}/2).
$$
Here $p^{B(x_0, \sqrt{t})}$ is the transition density function for the part process
$X^{B(x_0, \sqrt{t})}$ of $X$ killed upon leaving the ball $B(x_0, \sqrt{t})$.

\item{(v)} (Lemma \ref{L:compa}): Let $R\leq 1$ and $\delta <1$.
$Q_1=[t_0+2\delta R^2/3, \, t_0+ 5\delta R^2]\times B(x_0,
3R/2)$, $Q_2=[t_0+ \delta R^2/3, \, t_0+11\delta R^2/2]\times
B(x_0, 2R)$ and define
$Q_-$ and $Q_+$ as in Theorem \ref{T:PHI}. Let $h: [0,\infty)\times
\R^d\to \R_+$ be bounded and supported in $[0,\infty)\times
B(x_0,3R)^c$. Then there exists $c_5=c_5(\delta)>0$ such that
\[\E^{(t_1,y_1)}[h(  Z_{\tau_{Q_1}})]\le c_5\E^{(t_2,y_2)}
[h(  Z_{\tau_{Q_2}})] \qquad \mbox{for }  (t_1,y_1)\in Q_- \mbox{ and
}  (t_2,y_2)\in Q_+.\]
\end{description}

The proof of (iv)
uses ideas from \cite{BBCK}, where a similar
inequality is established for finite range pure jump process. However,
some difficulties arise due to the presence of the diffusion part.

The upper bound heat kernel estimate in Theorem \ref{mainHK} is established
by using method of scaling, by Meyer's construction of the process $X$ based
on finite range process $X^{(\lambda)}$, where the jumping kernel $J$ is replaced
by $J(x, y)\1_{\{|x-y|\leq \lambda\}}$, and by Davies' method from \cite{CKS}
to derive
an upper bound estimate for the transition density function of $X^{(\lambda)}$
through carefully chosen   testing functions. Here we need to select the value
of $\lambda$ in a very careful way that depends on the values of $t$ and $|x-y|$.

To get the lower bound heat kernel estimate in Theorem \ref{mainHK}, we
need a full scale parabolic Harnack principle
that extends Theorem \ref{T:PHI}
to all $R>0$ with the scale function $\tp (R):=R^2 \wedge \phi (R)$ in place of
$R\mapsto R^2$ there.  To establish such a full scale parabolic Harnack principle,
we show the following.

 \begin{description}
 \item{(iii')} Strengthened version of (iii) (Lemma \ref{L:4.12}):
  for every $x\in \R^d$, $r>0$ and any compact subset $A\subset
  Q(0,x,r) :=[0,\gamma_0 \tp(r)]\times B(x,r)$,
 $$    \bP^{(\gamma_0 \tp (r), x)}
 (\sigma_A <\tau_r) \geq c_3
 \frac{m_{d+1}(A)}{r^d \tp (r)}.
 $$
 Here $\gamma_0$ denotes the constant $\gamma (1/2, 1/2)$ in Proposition \ref{P:4.2}.

 \item{(vi)}
 (Corollary \ref{C:4.13}):
 For
 every $\delta \in (0, \gamma_0]$,
 there is a constant $c_6=c_6(\gamma)$
 so that for every $0<R\leq 1$, $r\in (0, R/4]$
 and $(t, x)\in Q(0, z, R/3)$ with
 $0<t\leq \gamma_0 \tp (R/3)-\delta \tp (r)$,
 $$
  \bP^{(\gamma_0\tp (R/3), z)}
  (\sigma_{U(t, x, r)}<\tau_{Q(0, z, R)})
 \geq c_6 \frac{r^d \tp (r)}{R^d \tp (R)},
 $$
 where $U(t, x, r):=\{t\}\times B(x, r)$.
\end{description}
With the full scale parabolic Harnack inequality, the lower bound heat kernel estimate
can then be derived once the following estimate is obtained.

\begin{description} \item{(vii)}
Tightness result (Proposition \ref{lower2}): there are constants
$c_7\geq 2$ and  $c_8>0$
such that for every $t>0$ and $x, y\in \R^d$ with
$|x-y|\geq c_7 \tp (t)$,
$$
\bP_x \left(X_t\in B(y, c_7\tp^{-1}(t)) \right) \geq c_8
\frac{t (\tp^{-1}(t))^d}{|x-y|^d \tp (|x-y|)}.
$$
\end{description}

Throughout the paper, we will define and use various Dirichlet forms, the corresponding processes
and heat kernels. For the convenience of the reader, we list the notations here.

\medskip

\halign{ #\hfil  \quad &#\hfil \quad &#\hfil \quad &#\hfil \cr
(Heat kernel)& (Process) & (Jump kernel) & (Dirichlet form)\cr
$p(t,x,y)$ &  $X$ & $J(x,y)$ &  $(\sE,\sF)=(\sE, W^{1,2}(\R^d))$\cr
$p^B(t,x,y)$ &  $X^B$  & $J(x,y)$ &  $(\sE,\sF^B)$: $X$ killed on exiting $B$\cr
$p^{(\lambda)}(t,x,y)$ &  $X^{(\lambda)}$ & $J(x, y) \1_{\{|x-y|\leq \lam\}}$ &
 $(\sE^{(\lam)}, W^{1,2}(\R^d))$\cr
 $p^{(\lambda; n)}(t,x,y)$ &  $X^{(\lambda; n)}$ &
 $J(x, y) \1_{\{|x-y|\leq \lam\}}\1_{B(n)\times B(n)}$
 & $(\sE^{(\lam; n)}, {\cal F}^{(\lam; n)})$\cr
$p_Y(t,x,y)$& $Y$ & $\kappa (x,y)|x-y|^{-d-\beta}$& subordinated Dirichlet form
 $\dashrightarrow
  (A)$\cr
$q^\delta(t,x,y)$ &  $Z^\delta$ & $J_\delta(x,y)
 \dashrightarrow
 (B)$ &  $(\sE^\delta,\sF^\delta)$\cr
$q^{\delta, B_r}(t,x,y)$ &  $Z^{\delta, B_r}$ & $J_\delta(x,y)$ &  $(\sE^\delta,\sF^{\delta, B_r})$:
$Z^\delta$ killed on exiting $B_r$\cr
$q_r^{\delta,B}(t, x, y)$ &  $r^{-1}Z^{\delta,B_r}_{r^2\,\cdot}$ & $J_\delta^{\angel{r}}(x,y)
 \dashrightarrow
(C)$ &  $(\sE^{\angel{r}}, \sF^{\angel{r},B})$:
$r^{-1}Z^\delta_{r^2\,\cdot}$ killed on exiting $B$\cr
$p_r(t, x, y)$ &  $X^{\angel{r}}$ & $J^{\angel{r}}(x,y)
 \dashrightarrow
(D)$ &
$(\sE^{\angel{r}}, \sF^{\angel{r}})=(\sE^{\angel{r}},W^{1,2}(\R^d))$\cr
$p^{(\lambda)}_r(t, x, y)$ &  $X^{\angel{r,\lambda}}$ & $J^{\angel{r}}(x,y)\1_{\{|x-y|\le \lambda\}}$ &
$(\sE^{\angel{r,\lambda}}, W^{1,2}(\R^d))$ \cr
}

\medskip
where in the above,
\begin{description}
\item{(A)}
$Y$ is the subordination of the symmetric diffusion for $\nabla (A\nabla)$, the local part of
$\sE$, by the subordinator $\eta=\{t+c_0\eta_t^{(1)}, t\ge 0\}$, where
$\{\eta_t^{(1)}\}$ is a $(\beta/2)$-subordinator.

\item{(B)}
$J_\delta(x,y):=J(x,y) \1_{\{|x-y|\ge \delta\}}
+\kappa (x,y)|x-y|^{-d-\beta} \1_{\{|x-y|< \delta\}}.$

 \item{(C)}
$q_r^{\delta,B}(t, x, y)=q_r^B(t, x, y):=r^dq^{\delta, B_r}(r^2t,rx,ry)$,
$Z^{\angel{r}}_t:=r^{-1}Z^\delta_{r^2t}$,
$J^{\angel{r}}_\delta(x,y):=r^{d+2}J_\delta(rx,ry)$ for $r\in (0,1]$.

\item{(D)}
$p_r(t, x, y):=r^{d}p(\tp(r)t, rx, ry)$,
$X_t^{\angel{r}}:=r^{-1}X_{\tp(r)t},~J^{\angel{r}}(x,y):=\tp(r)r^dJ(rx,ry)$
 for $r>0$.
\end{description}

\section{Heat kernel upper bound estimate and exit time estimate}

Throughout this paper, We always assume the uniform elliptic condition \eqref{unielli}
holds for the diffusion matrix $A$.
Let $(\sE, \sF)$ be the Dirichlet form in \eqref{e:DF} with the
jumping kernel $J$ satisfying the conditions
\eqref{e:J1} and \eqref{e:J2}. We start this section by giving a

\medskip

\noindent{\bf Proof of Proposition \ref{P:charaF}}:
 For any
$u\in C^1_0 (\R^d)$, we have
\[\int_{\R^d}
\nabla u(x)\cdot A(x) \nabla u(x)dx+\|u\|_2^2\asymp \int_{\R^d} |\nabla u(x)|^2dx+\|u\|_2^2=:{\cal C}_{1,c}(u,u),\] and
\begin{eqnarray}
\int_{\mathbb R^d}(u(x)-u(y))^2J(x,y)dxdy&\le& \int_{|x-y|\le 1}
(u(x)-u(y))^2J(x,y)dxdy+c_1\|u\|_2^2\nonumber\\
&\le& c_2\Big(\int_{\mathbb R^d}
\frac{(u(x)-u(y))^2}{|x-y|^{d+\beta }}dxdy+\|u\|_2^2\Big)
=:c_2{\cal C}_{1,d}(u,u).~~~~~~~\label{eq:1_6nij}\end{eqnarray}
Using Fourier transform, it is well-known that
 \begin{equation}\label{e:Fourhat}
  {\cal C}_{1,d}(u,u) = c \int_{\R^d} (|\xi|^\beta+1) |\wh u (\xi)|^2
d\xi \leq 2 c \int_{\R^d} (|\xi|^2+1) |\wh u (\xi)|^2 d\xi
 = c_3 {\cal C}_{1,c}(u,u).
 \end{equation}
Thus we have $\sE(u,u)\asymp {\cal C}_{1,c}(u,u)$ for all $u\in
C^1_0 (\R^d)$. It follows then
\[\sF=\overline{C^1_0 (\R^d)}^{\sE_1}=
\overline{C^1_0 (\R^d)}^{{\cal C}_{1,c}}=W^{1,2}(\R^d).\]
 \qed

\subsection{Heat kernel upper bound estimate}

By the Nash's inequality
 \begin{equation}\label{e:Nash}
\|f\|_2^{2+4/d}\le c_1\int_{\mathbb R^d} |\nabla
u(x)|^2dx\cdot\|f\|_1^{4/d}\le c_2\sE (f,f)\|f\|_1^{4/d}\qquad
\hbox{for }  f\in W^{1,2}(\R^d),
 \end{equation}
 we have, by Theorem
\cite[Theorem 2.9]{CKS} and \cite[Theorem 3.1]{BBCK}, that there is
a properly $\sE$-exceptional set ${\cal N}\subset \R^d$ of $X$ and
a positive symmetric kernel $p(t, x, y)$ defined on $[0, \infty)
\times (\R^d\setminus {\cal N})\times (\R^d\setminus {\cal N})$
such that for every $x\in \R^d\setminus {\cal N}$ and $t>0$,
$$ \E_x \left[ f(X_t) \right] =\int_{\R^d} p(t, x, y) f(y) m_d (dy),
$$
$$ p(t+s, x, y)=\int_{\R^d} p(t, x, z)p(s, z, y) \qquad \hbox{for
every } t, s>0 \hbox{ and } x, y \in \R^d \setminus {\cal N},
  $$
and
 \begin{equation}\label{e:p1}
 p(t, x, y) \leq c t^{-d/2} \qquad \hbox{for } t>0 \hbox{ and every }
   x, y \in \R^d \setminus {\cal N}.
 \end{equation}
Moreover, there is an $\sE$-nest $\{F_k, k\geq 1\}$ of compact
subsets of $\R^d$ so that $\displaystyle {\cal N}= \R^d \setminus
\cup_{k=1}^\infty F_k$ and that for every $t>0$ and $y\in \R^d
\setminus {\cal N}$, $x\mapsto p(t, x, y)$ is continuous on each
$F_k$. Later, as a consequence of the H\"older continuity result for
parabolic functions, $p(t, x, y)$ in fact has a continuous version
so the exceptional set ${\cal N}$ can be taken to be an empty set.

\medskip

Now, for $\lam\in \bQ_+$, where $\bQ_+$ is the set of positive rational
numbers, let $(\sE^{(\lam)}, W^{1,2}(\R^d))$ be the Dirichlet form defined
by \eqref{e:DF} but with the jumping kernel $J(x, y) \1_{\{|x-y|\leq
\lam\}}$   in place of $J(x, y)$. Let $X^{(\lam)}$ be the symmetric
strong Markov process associated with $(\sE^{(\lam)}, W^{1,2}(\R^d))$,
and let $p^{(\lam)} (t, x, y)$ be its transition density function.
\begin{propn}\label{thm:newtkjuly}
Let $\displaystyle{ \delta (\lam):= \sup_{\xi \in \R^d}
\int_{\{\eta \in \R^d: \, |\eta-\xi|\leq \lam\}} |\xi -\eta|^2 J
(\eta,\xi) d \eta}$. Then, there exist $c_1,c_2>0$ (independent of $\lam\in \bQ_+$)
such that for any $s>0$, the following holds for all $t>0$ and q.e. $x,y$,
\begin{equation}\label{qdnief}
 p^{(\lam)}(t, x, y) \leq c_1 t^{-d/2} \exp \left( - s |x-y|
  + c_2 s^2 \left(1+ e^{2\lam s} \delta (\lam) \right) t \right).
\end{equation}
\end{propn}
\pf
First, note that by condition \eqref{e:J1}, we have
\begin{equation}\label{e:delta}
 \lim_{\lam\to 0} \delta (\lam) =0.
 \end{equation}
We use Davies' method to derive the desired heat kernel upper bound.
From Nash's inequality
\eqref{e:Nash}, by the same reasoning as that for $X$ at the
beginning of this section, the symmetric process $X^{(\lam)}$ has a
quasi-continuous transition density function $p^{(\lam)} (t, x, y)$ defined
on $[0, \infty)\times (\R^d\setminus \NN_\lam ) \times ( \R^d \setminus
\NN_\lam)$ such that
\begin{equation}\label{e:q_twe}
 p^{(\lam)}(t, x, y) \leq c_1 \, t^{-d/2} \qquad \hbox{for every } t>0 \hbox{ and }
 x, y \in  \R^d\setminus \NN_\lam.
 \end{equation}
 Note that the above constant $c_1>0$ is independent of $\lam>0$.
 By (\ref{e:Fourhat}),
 we have $\sE^{(\lam)}_1(u,u)\asymp {\cal C}_{1,c}(u,u)\asymp
\sE_1(u,u)$, so a set is $\sE^{(\lam)}_1$-exceptional if and only if it is
 $\sE_1$-exceptional. Thus, letting $\NN=\cup_{\lam\in\bQ_+}\NN_\lam$,
$\NN$ is a $\sE_1$-exceptional set.
(\ref{e:q_twe}) together with \cite[Theorem 3.25]{CKS} and
 \cite[Theorem 3.2]{BBCK}   implies that there exist constants $C>0$ and $c>0$,
 such that
\begin{equation}\label{e:q}
p^{(\lam)} (t,x,y)\leq c_1\, t^{-d/2} \;\exp\left( -|\psi(y)-\psi(x)|+C\;
\Lambda_\lam (\psi)^2\; t\right)
\end{equation}
for all $t>0$,  $x,y\in \R^d\setminus \NN$, and for any function $\psi$ having
 $\Lambda_\lam (\psi)<\infty$.
Here
\begin{equation*}
\Lambda_\lam (\psi)^2= \|e^{-2\psi}\Gamma_{\lam}[e^\psi]\|_\infty \vee
\|e^{2\psi}\Gamma_{\lam}[e^{-\psi}]\|_\infty .
\end{equation*}
where  for $\xi \in \R^d$,
\begin{equation}  \label{e:densi}
\Gamma_\lam [v] (\xi) := \sum_{i,j=1}^da_{ij}( \xi)\frac{\partial
v}{\partial x_i} (\xi)\frac{\partial v}{\partial x_j}(\xi)+
\int_{\{\eta \in \R^d: \, |\eta-\xi|\leq \lam\}} (v(\eta)-v(\xi))^2J
(\eta,\xi) d \eta,
\end{equation}
 For $s>0$, take
$$
\psi(\xi):=s\,  \left(  |\xi-x| \wedge |x- y| \right) \qquad
\hbox{for } \xi \in \R^d.
$$
Note that
$|\psi(\eta)-\psi(\xi)|\leq s\,
 |\eta-\xi|$ for all $\xi,\eta\in \R^d$. So for $\xi \in \R^d$,
\begin{align*}
e^{-2\psi(\xi)}\Gamma_\lam [e^\psi](\xi)&\le c_2|\nabla \psi (\xi)|^2+
\int_{|\eta-\xi|\leq \lam} (1-e^{\psi(\eta)-\psi(\xi)})^2J (\eta,\xi)
d \eta\\
&\leq c_2s^2 +\int_{|\eta-\xi|\leq \lam} (\psi(\eta)-\psi(\xi))^2\;
e^{2|\psi(\eta)-\psi(\xi)|}J (\eta,\xi) d \eta\\
&\leq c_2s^2 + s^2\, e^{2\lam s  } \int_{|\eta-\xi|\leq
\lam}|\eta-\xi|^2 J (\eta,\xi) d\eta\\
&\leq c_2 s^2 \left(1+ e^{2\lam s }\delta (\lam) \right) .
\end{align*}
Here $c_2>0$ is independent of $\lam\in \bQ_+$. The same estimate holds for
$e^{2\psi(\xi)}\Gamma_\lam [e^{-\psi}](\xi)$.
 So we have the desired estimate. \qed

\subsection{Conservativeness}

\begin{thm}\label{consernew}
The process $X$ is conservative; that is, $X$ has infinite lifetime.
\end{thm}

\pf Recall the process $X^{(\lam)}$ defined in the previous subsection. $X$ can be obtained from
$X^{(\lam)}$ through Meyer's construction by adding all the jumps
whose size is larger than $\lam$ (see Remarks 3.4-3.5 of \cite{BBCK}
and Lemma 3.1 of \cite{BGK}).
Note that by \eqref{e:J1} and \eqref{e:J2}, there is a   constant
$b_0>0$ such that
\begin{equation}\label{e:Je1}
\sup_{x\in \R^d} \int_{\R^d} \1_{\{|x-y|>
\lam\}} J(x, y) dy \leq b_0  \lam^{-\beta}
\qquad \hbox{for every } \lam\in (0, 1].
\end{equation}
Thus, it suffices to show that $X^{(\lam)}$ is conservative.
To show this, we look at reflected
 jump-diffusions
with jumping kernel $J(x, y) \1_{\{|x-y|\leq
\lam\}}$ in big balls, as in \cite[Theorem 4.7]{CK2}. In the following, we fix $\lam
\in \bQ_+$.
Let $x_0\in \R^d$, $r_n\geq 100\lam$. Define $B(n)=\overline{B(x_0, r_n)}$ and
\begin{eqnarray*}
\sE^{(\lam; n)}(f,f)&=&\int_{B(n)} \nabla f(x)\cdot A(x) \nabla
f(x)dx+\int_{B(n)}\int_{B(n)}(f(x)-f(y))^2
J(x, y) \1_{\{|x-y|\leq \lam\}} dx dy,\\
\sF^{(\lam; n)}&=&\overline{\{
 f\in C^1 (\overline {B(n)}): \
\sE^{\ (\lam; n)}(f,f)<\infty\}}^{\sE^{(\lam; n)}_1},
\end{eqnarray*}
where $\sE^{(\lam; n)}_1(u, u):=\sE^{(\lam; n)}(u, u)+\int_{B(n)} u(x)^2dx$.
 Clearly  $(\sE^{(\lam; n)}, {\cal F}^{(\lam; n)})$ is a regular symmetric Dirichlet form
on $L^2( \overline{B(n)}; dx)$.
Let $X^{(\lam; n)}$ be the Hunt process on $\overline {B(n)}$
  associated with $(\sE^{(\lam; n)},\sF^{(\lam; n)})$. Since a constant
function $1 \in {\cal F}^{(\lam; n)}$ with $\sE^{(\lam; n)} (1, 1)=0$, $X^{(\lam; n)}$
is recurrent and so $X^{(\lam; n)}$ is conservative. Let $p^{(\lam; n)} (t, x, y)$ be
the transition density function of $X^{(\lam; n)}$. Then, similarly to the
proof of Proposition \ref{thm:newtkjuly}, we see that $p^{(\lam; n)} (t, x, y)$ exists
for all $t>0$, $x,y\in B(n)\setminus {\cal N}_n$, where ${\cal N}_n$ is a properly exceptional set
for  $X^{(\lam; n)}$, and
 moreover
it enjoys the estimate (\ref{qdnief}) with constants independent of $n$.
Using (\ref{qdnief}) with $s=1$, for $x\in B(n)\setminus {\cal N}_n$, $t\in [1,2]$ and $R\le r_n$, we have
\begin{eqnarray*}
\bP_x \left( | X^{(\lam;n)}_{s}-x| \ge R\right)&=&\int_{B(n)\setminus B(x,R)}
p^{(\lam; n)} (t, x, y)dy\\
&\le & c_1\int_{B(n)\setminus B(x,R)}e^{-|x-y|}dy\le c_2e^{-R},
\end{eqnarray*}
where $c_1,c_2$ may depend on $\lam$, but they are independent of $n$ and $R$.
Given this estimate, the rest is the same as
that of \cite[Theorem 4.7]{CK2}. We will sketch the argument.
Note that
for $x\in B_{r_n-\lam}\setminus {\cal N}_n$, $X^{(\lam; n)}$ has the same distribution
as that of  $X^{(\lam)}$ before $X^{(\lam;n)}$ leaves the ball $B_{r_n-\lam}$. Thus, estimating
as in \cite[(4.23)]{CK2}, we have for a.e. $x\in B_{r_0}$,
\begin{eqnarray*}
\bP_x \left(\zeta>1 \mbox{ and } \sup_{s\le 1}| X^{(\lam)}_{s}-x| \le R\right)&
\ge&\bP_x \left(\sup_{s\le 1}| X^{(\lam;n)}_{s}-x| \le R\right)\\
&\ge & 1-2c_2e^{-R/2}\qquad\mbox{for every } R>0,
\end{eqnarray*}
where $\zeta$ is the lifetime of  $X^{(\lam)}$.
Passing $R\to \infty$, we have for a.e. $x\in B_{r_0}$,
\begin{equation}\label{e:dwhibqq}
\bP_x(X^{(\lam)}_1\in \R^d)=1.
\end{equation}
Taking $r_0\uparrow \infty$, (\ref{e:dwhibqq}) holds for  a.e. $x\in \R^d$; by the Markov property,
$\bP_x(X^{(\lam)}_t\in \R^d)=1$ for every rational $t>0$.
Since for each rational $t>0$, $P^{(r)}_t 1$ is finely continuous
and $P^{(r)}_t 1=1$ a.e. on $\R^d$, we must have $P^{(r)}_t 1=1$
 q.e. on $\R^d$, so that $\bP_x(\zeta =\infty)=1$ for q.e. $x\in \R^d$.  \qed

\subsection{Exit time estimate}

For $A\subset \R^d$, denote by
$$ \tau_A:=\inf\{t>0: X_t\notin A\}
$$
the first exit time from $A$ by $X$.

\begin{lem}\label{L:exit}
For every $x_0\in \R^d$ and $r>0$,
$\E_x \left[ \tau_{B(x_0,r)} \right] \leq c_1 r^2$ for every $x\in
B(x_0, r)\setminus {\cal N}$.
\end{lem}

\pf
The proof for this is nowadays standard, see for example \cite{Ch}.
For reader's convenience, we spell out the details here.
Let $c>0$ be the constant in  \eqref{e:p1}. Take $c_2>0$ be large enough so that
 $$c \, m_d (B(0, 1)) \, c_2^{-d/2}\leq \tfrac12.
 $$
 Then for every $r>0$, $x_0\in \R^d$ and $x\in B(x_0, r) \setminus {\cal N}$,
  with $t:=c_2r^2$ we have by \eqref{e:p1},
\[ \bP_x (X_t\in B(x_0,r))=\int_{B(x_0,r)} p(t,x,z)dz
\leq c\, t^{-d/2} m_d(B(x_0,r))  \leq \tfrac12. \]
Since $X$ is conservative,  this implies
that for every $x\in B(x_0, r)\setminus {\cal N}$,
$$ \bP_x ( \tau_{B(x_0,r)}\leq  t) \geq \bP_x (X_t \notin B(x_0, r))
\geq 1/2.
$$
In other words, we have $\bP_x(\tau_{B(x_0,r)}> t)\leq
\tfrac12$.
  By the Markov property of $X$, for integer $k\geq 1$,
\[ \bP_x (\tau_{B(x_0,r)}> (k+1)t)\leq
  \E_x[\bP_{X_{kt}}(\tau_{B(x_0,r)}> t);
 \tau_{B(x_0,r)}>mt]\leq \tfrac12 \bP_x(\tau_{B(x_0,r)}>kt). \]
Using mathematical induction, we can conclude that for every $k\geq
1$,
\[ \bP_x (\tau_{B(x_0,r)}>kt)\leq 2^{-k}, \]
which yields the desired estimate $\E_x \left[ \tau_{B(x_0,r)}
\right]
 \leq \sum_{k=0}^\infty t \bP_x(\tau_{B(x_0,r)}>kt)
 \leq c_1 r^2$.
\qed

\begin{lem}\label{L:exit2}
There is are  constants $a_0, r_0\in (0, 1)$ so that for every $x\in
\R^d\setminus {\cal N}$,
$$ \bP_x \left( \sup_{s\leq a_0 r^2} |X_s-X_0| \leq r \right) \geq
1/4 \qquad \hbox{for every } r\in (0, r_0] .
$$
Consequently, there exists a constant $a_1>0$ so that for every
$x\in \R^d \setminus {\cal N}$,
$$\E_x \left[ \tau_{B(x, r)} \right] \geq a_1 r^2 \qquad \hbox{for
every } r\in (0, r_0].
$$
\end{lem}

\pf
By Lemma 3.6 of \cite{BBCK} and (\ref{e:Je1}),
we have for $0<r\leq 1$,
\begin{eqnarray*}
 \bP_x \left( \sup_{s\leq a_0 r^2} |X_s-X_0| \leq r \right)
  &\geq&  e^{- (b_0 r^{-\beta})( a_0 r^2) } \,
  \bP_x \left( \sup_{s\leq a_0 r^2} |X^{(r)}_s-X^{(r)}_0| \leq r \right) \\
  &\geq & e^{- a_0 b_0  }\,
   \bP_x \left( \sup_{s\leq a_0 r^2} |X^{(r)}_s-X^{(r)}_0| \leq r \right).
   \end{eqnarray*}
   So it suffices to show that there is a positive constant $a_0\in (0, 1)$ small
   so that $\leq \log 2$,
   \begin{equation}\label{e:a0}
     a_0b_0  < b_0 a_0^{\beta/2}<\log (8/7)
   \end{equation}
   and  that
   $$  \bP_x \left( \sup_{s\leq a_0 r^2} |X^{(r)}_s-X^{(r)}_0| \leq r \right) \geq 1/2
   \qquad \hbox{for every }  r\in (0, r_0]\cap \bQ
   \hbox{ and } x \in \R^d\setminus \NN.
   $$
Taking
$ s= 1/\sqrt{t}$ in (\ref{qdnief}), we have
   \begin{equation}\label{e:cal1}
    p^{(r)} (t, x, y) \leq c_0 t^{-d/2} \exp
\left( -   \frac {|x-y|}{\sqrt{t}}
  + c_2 \left(1+ e^{2 r/\sqrt{t}} \delta (r) \right)  \right).
  \end{equation}
Using polar coordinate,
 \begin{equation}\label{e:cal2}
  \int_{\{|x-y
|\geq r/2\}} c_0 t^{-d/2} e^{2c_2} \exp \left( -   \frac {
|x-y|}{\sqrt{t}} \right) dy = \omega_d c_0 e^{2c_1}
\int_{\frac{r}{2\sqrt{t}}}^\infty e^{-v} dv,
  \end{equation}
where $\omega_d$ is a positive constant that depends only on
dimension $d$. Let $a_0>0$ be small enough so that
$$ \omega_d c_0 e^{2c_2}
\int_{1/(2\sqrt{a_0})}^\infty e^{-v} dv < 1/8.
$$
Due to \eqref{e:delta}, there exists $r_0 \in (0, 1)$ so that
 $$ e^{2 /\sqrt{a_0}} \delta (r) \leq 1 \qquad \hbox{for every }
 r\in (0, r_0].
 $$
This together with \eqref{e:cal1} and \eqref{e:cal2} implies that
 for every $r\in (0, r_0]\cap \bQ$ and $ x\in \R^d$,
$$ \bP_x \left( | X^{(r)}_{a_0 r^2}-X^{(r)}_0| \geq r/2 \right)
  = \int_{\{|y-x|\geq r/2\}} p^{(r)} (a_0 r^2, x, y) dy \leq 1/8.
$$
Moreover, by \cite[Lemma 3.6]{BBCK}, we have for every $s\leq a_0
r^2$ with $r\in (0, r_0]\cap \bQ$,
\begin{eqnarray*}
 \bP_x \left( | X^{(r)}_{s}-x| < r /2 \right)
&\geq  & \bP_x \left( | X^{(r)}_{s}-x| < \sqrt{s/a_0}/2 \right) \\
&  \geq &   e^{-s  \, J_{s, r}} \, \bP_x \left( \big|
X^{(\sqrt{s/a_0})}_{s}-x \big|
 < \sqrt{s/a_0}/2 \right) \\
 &\geq &  \frac78 \,  e^{-s  \, J_{s, r}} ,
\end{eqnarray*}
where
$$ J_{s, r} =\sup_{x\in \R^d} \int_{\R^d} \1_{\{ \sqrt{s/a_0} < |x-y| \leq r\}}
J(x, y) dy.
$$
By \eqref{e:Je1} and \eqref{e:a0},
$$ s J_{s, r} \leq b_0 a_0^{\beta/2} s^{(2-\beta)/2} \leq b_0
a_0^{\beta/2} < \log (8/7)
$$
and so
$$ \inf_{x\in \R^d \setminus \NN} \bP_x \left( | X^{(r)}_{s}-x| < r /2 \right) \geq
  (7/8)^2 > 3/4.
$$
In other words, we have
$$ \sup_{x\in \R^d \setminus \NN }  \bP_x \left( | X^{(r)}_{s}-x| \geq  r /2 \right) < 1/4
\qquad \hbox{for every }  s\leq  a_0 r^2.
 $$
Now, since $X^{(r)}$ is conservative, by  Lemma 3.8 of \cite{BBCK},
$$ \sup_{x\in \R^d \setminus \NN }
 \bP_x \left( \sup_{s\leq a_0r^2} | X^{(r)}_s - X^{(r)}_0| \geq  r\right)
<1/2,
$$
for every $r\in (0, r_0]\cap \bQ$.
This proves the lemma. \qed

\section{Short time near-diagonal heat kernel
lower bound estimate}

Let $X$ be the strong Markov process associated with the Dirichlet
form $(\sE, \sF)$ of \eqref{e:DF} with the jumping kernel satisfying
the condition \eqref{e:J1}-\eqref{e:J2} and \eqref{e:J1b}.
Recall that $p(t, x, y)$ is the transition density function for $X$.
For a ball $B\subset \R^d$, denote by $p^B(t, x, y)$ the transition
density function of the subprocess $X^B$ of $X$ killed upon exiting $B$.
In this section we will establish the following.

\begin{thm}\label{T:3.4}
 For each $t_0>0$, there exists $c =c (t_0)>0$
 such that for every $x_0\in \R^d$ and $t\le t_0$,
  $$
p^{ B(x_0, \sqrt{t })} (t, x,y)\ge c\, t^{-d/2}\qquad \hbox{for q.e. } x,
y\in B(x_0, \sqrt{t }/2)
 $$
 and
\[ p (t, x,y)\ge c \, t^{-d/2}\qquad
\mbox{for q.e. $x, y$ with } |x-y|^2 \leq t.
 \]
\end{thm}

This result will be used in later sections   with $t_0=1$.
For its proof, we adopt an approach from \cite{BBCK} that deals
with finite range pure jump processes.
But there are some new technical difficulties to overcome
in our setting.

Fix $x_0\in \R^d$ and let $a_1=12/(2-\beta )$. (In fact,
the following argument works for any fixed $a_1$ bigger than
$4\vee (6/(2-\beta ))$.)  For  $r>0$, define
\[\Psi_r(x)=c((1-r^{-1}|x-x_0|)_+)^{a_1},\]
where $c>0$ is the normalizing constant such that $\int_{\mathbb
R^d}\Psi_r(x)dx=1$. Then the following weighted Poincar\'e
inequality holds. (See, for example, \cite[Theorem 5.3.4]{SC} for
the proof.)
\begin{propn}\label{wtedPI} There is a positive constant
$c_1=c_1(d)$
independent of $r$, such that
\[ \int_{B(x_0,r)} ( u(x)-u_{\Psi_r})^2 \Psi_r (x) dx \\
\leq c_1 r^2\int_{B(x_0,r)} |\nabla u(x)|^2\Psi_r (x) dx
\qquad \hbox{for }
u\in C^\infty_b(\R^d).\]
Here $u_{\Psi_r} := \int_{B(x_0,r)} u(x) \Psi_r (x) dx$.
\end{propn}

\bigskip

  Let $W$ be the symmetric diffusion
that corresponds to the divergence form operator $\nabla (A\nabla
)$, the local part  of $\sE$. Let $\eta^{(1)}=\{ \eta^{(1)}_t, t\geq
0\}$ be an $(\beta/2)$-subordinator and define $\eta_t=t+ c_0
\eta^{(1)}_t$, where $c_0>0$ is a large constant to be chosen at the
end of this paragraph. Define $ Y$ to be the subordination of $W$ by
the subordinator $\eta=\{ \eta_t; t\geq 0\}$.   Note that $ Y$ is a
symmetric strong Markov process, whose continuous part has the same
law as $W$, and its jumping part comes from the subordination of $W$
by $c_0 \eta^{(1)}$. By the uniform ellipticity (\ref{unielli}) of
the diffusion matrix $A(x)$, the heat kernel of $W$ enjoys
Aronson-type two-sided Gaussian estimate. It follows that  (see
\cite{Stos}) the jump kernel of $ Y$ is of the form $\kappa
(x,y)/|x-y|^{d+\beta}$, where $\kappa (x, y)$ is a symmetric
measurable function that is bounded between two positive constants.
By taking $c_0>0$   sufficiently large, we can and do assume that
  $$ J(x, y) \leq \frac{\kappa (x, y)}{ |x-y|^{d+\beta}}
   \qquad \hbox{for all } |x-y| \leq 1.
  $$

For $\delta \in (0, 1)$, set
\begin{equation}\label{eqn:3.7}  J_\delta (x, y)=\
\begin{cases}
J(x, y) &\hbox{for } |x-y|\geq \delta; \\
 \kappa (x,y) |y-x|^{-d-\beta} &\hbox{for }  |x-y|<\delta ,
\end{cases}
\end{equation}
and define $(\sE^\delta, \sF^\delta)$ with $J_\delta$ in place of
$J$ in the definition of $(\sE, \sF)$.

For $\delta \in (0, 1)$, let $Z^\delta$ be the symmetric Markov
process associated with $(\sE^\delta, \sF^\delta)$. Note that the
jumping kernel for $Z^\delta$ differs from that of $ Y$ by a bounded
and integrable kernel. So $Z^\delta$ can be constructed from $  Y$
through Meyer's construction (see Remarks 3.4 and 3.5 of
\cite{BBCK} and Lemma 3.1 of \cite{BGK}). Consequently,
the process $Z^\delta$ can be modified
to start from every point in $\R^d$ and $Z^\delta$ is conservative.
Moreover by a similar proof to that in \cite{BBCK}, we can show that
 $Z^\delta$ has a quasi-continuous transition density function
$q^\delta (t, x, y)$ defined on $[0, \infty) \times \R^d \times
\R^d$, with respect to the Lebesgue measure on $\R^d$. Since $  Y$
is a subordination of $W$, we can readily get a two-sided kernel
estimate on $p_Y (t, x, y)$ of $Y$ from that of $W$. In fact, since
the heat kernel of $W$ is comparable to that of Brownian motion,
$p_Y(t, x, y)$ is comparable to that of the independent sum of
Brownian motion and a rotationally  symmetric $\beta$-stable
process. So by \cite{SV07},
\begin{eqnarray}
 && c_1 \left( t^{-d/2}\wedge t^{-d/\beta} \right)  \left( t^{-d/2} e^{-c_2|x-y|^2/t}
+ t^{-d/\beta} \left( 1\wedge \frac{t}{|x-y|^{d+\beta}} \right) \right) \label{e:est2}\\
&\leq & p_Y (t, x, y) \leq   c_3 \left( t^{-d/2}\wedge
t^{-d/\beta} \right) \left( t^{-d/2} e^{-4_2|x-y|^2/t} +
t^{-d/\beta} \left( 1\wedge \frac{t}{|x-y|^{d+\beta }} \right)
\right)  \nonumber
\end{eqnarray}
for all $t>0$ and $x, y\in \R^d$. Consequently, parabolic Harnack
principle holds for $Y$ (see \cite[Theorem 4.5]{SV07}). On the other
hand, as a consequence of Meyer's construction (see the proof of
Proposition 2.1 of \cite{CKK}) and \eqref{e:est2}, there are
constant $t_0, r \in (0, 1)$ and $c>1$, which depend on $\delta$, so
that
 \begin{equation}\label{e:est1}
  c^{-1} p_Y (t, x, y) \leq q^\delta (t,
x, y) \leq c \, p_Y(t, x, y) \qquad \hbox{for } t\in (0, t_0]
\hbox{ and } |x-y| \leq r_0.
  \end{equation}
  From \eqref{e:est1}, we can easily show that parabolic Harnack principle
  holds at small-size scale
  for $Z^\delta$ and that its parabolic functions are jointly
  continuous (see \cite[Remark 4.3(ii)]{CKK}).
   In particular,
  $q^\delta (t, x, y)$ is jointly continuous on $\R_+ \times \R^d
  \times \R^d$.

\medskip

For $r\in (0,1]$, let $B_r=B(0,r)$ and let $(\sE^\delta,
\sF^{\delta,B_r})$ be the Dirichlet form corresponding to the
process $Z^\delta$ killed on leaving the ball $B_r$. Let
$q^{\delta,B_r}(t,x,y)$ be its heat kernel with respect to the
Lebesgue measure in $B_r$. We first prove the following, which
corresponds to Lemmas 4.5, 4.6 and 4.7 in \cite{BBCK}.

\begin{propn}\label{thm:bbcklem47}
(i)  For each $t>0$ and $y_0\in B_r$, we have
\[q^{\delta,B_r}(t,\cdot,y_0), \ \
\frac{\Psi_r (\cdot )}{q^{\delta,B_r}(t, \cdot, y_0)}\in \sF^{\delta,B_r}.\]
(ii) Fix $y_0\in B$ and
let $\bar G(t)=\int_{B_r} \Psi_r(x) \log q^{\delta,B_r}(t, x, y_0) \, dx.$
Then for every $t>0$,
\[\bar G'(t)=-\sE \left(q^{\delta,B_r}(t, \cdot , y_0), \,
\frac{\Psi_r (\cdot )}{ q^{\delta,B_r}(t, \cdot, y_0)} \right).\]
\end{propn}

\medskip

The following lemma plays a key role in our proof of above
proposition.

\begin{lem}\label{thm:bbckpro43}
Assume $0<\delta <1/16$.  Let  $0<t_1<t_2 < \infty$ and $r\in
(16\delta , 1]$. There is a constant $c_1=c_1(\delta, r, t_0,t_1)>0$
such that
$$ q^{\delta,B_r}(t, x, y)\geq c_1(r-|x|)^2 (r-|y|)^2 \qquad
\mbox{ for every $ t\in [t_1,t_2]$ and  $x, y \in B_r$.} $$
\end{lem}

\pf
Due to the Chapman-Kolmogorov equation, without loss of
generality, we can and do assume that
$$t_1< 3 a_0  \min \{\delta_0 r, r_0\}^2/16, $$
where $\delta_0 \in (0, 1)$ is the constant in
\eqref{e:J1} and \eqref{e:J1b}.
and $a_0$ and $r_0$ are the constant in Lemma
\ref{L:exit2}.

First, since as mentioned above $Z^\delta$ enjoys parabolic Harnack
principle at the small-size scale, we have by the same proof as that
for Lemma 4.2 of \cite{BBCK} that for every $\gamma \in (0, 1)$,
there is a constant $c_\gamma >0$ so that
  \begin{equation}\label{e:ql}
   q^{\delta, B_r} (t, x, y) \geq c_\gamma \qquad \hbox{for } t\in [
   t_1/12, \,  t_2  ] \hbox{ and } x, y \in B(0, \gamma r).
 \end{equation}
 So it suffices to prove the lemma for $x, y \in B_r$ with
 $$ \max\{r-|x|, r-|y|\} < r_1:=\min\{r_0, \, \delta_0 r/8, \,  t_1/(4a_0)\}.
 $$

 Let $y\in B_r$ with
$\delta (y):=r-|y|<r_1 $. Take $y_0\in B(0, (1-3\delta_0/4 )r)$ with $|y-y_0|=
 \delta_0 r$. Define
   $T:=\inf\{t>0: \ |Z^\delta_t-Z^\delta_{t-}|\geq \delta_0r \}$
  and set  $s_0=t_1/3$. By the strong Markov
  property of $Z^\delta$,

\begin{eqnarray}
 && \bP_y \left( Z^\delta_{s_0} \in B(0, (1-\delta_0/2)r)
 \hbox{ and }  \tau_{B_r} > s_0)\right) \nonumber \\
&\geq &  \bP_y \left(T\leq a_0 \delta (y)^2/4,
 \  Z^{\delta}_T    \in B(y_0 , \, \delta_0 r/16), \ \sup_{s<T}
  |Z^\delta_{s}-y|\leq  \delta (y)/2  \right. \nonumber \\
 && \hskip 0.3truein \left.
   \hbox{ and } \sup_{s\in [T, s_0+T]} |Z^\delta_s-Z^\delta_T| \leq \delta_0 r/4 \right) \nonumber  \\
 &\geq & \bP_y \left(T\leq a_0 \delta (y)^2/4,
 \  Z^{\delta}_T   \in B(y_0 , \, \delta_0 r/16) \hbox{ and } \sup_{s<T}
 |Z^\delta_{s}-y|\leq
 \delta (y)/2  \right) \nonumber \\
    && \cdot \inf_{y\in \R^d \setminus \NN} \bP_x \left(
    \sup_{s\in [0, s_0]} |Z^\delta_s-x| \leq \delta_0 r/4 \right)  \label{e:pest}.
    \end{eqnarray}
Note that by conditions \eqref{e:J1}-\eqref{e:J2} and \eqref{e:J1b},
 $$ \kappa_1
:=\sup_{x\in \R^d} \int_{\R^d} \1_{\{  |x-z|> \delta_0r \}}
J_\delta (x, z) dz  <\infty
$$
and
 $$ \kappa_2 :=
\inf_{y\in B_r} \, \inf_{x\in B(y, \, \delta_0 r/16)\setminus \NN }\int_{B(y_0 , \, \delta_0 r/16)}
  J_\delta (x, z) dz >0.
  $$
  As $T$ is the first time the process $Z^\delta$ makes a jump of size no less than
  $\delta_0r$, $T$ is stochastically dominated from above
  by the exponential random variable with parameter $\kappa_1$
  and at time $T$, process $Z^\delta$ jumps
  to position $z$ according
  to the probability kernel
  $$ \frac{J_\delta (Z^\delta_{T-}, dz)}{
   \int_{\{w: \, |w-Z^\delta_{T-}| \geq \delta_0 r) \}} J_\delta (Z^\delta_{T-}, dw) }.
  $$
Thus we have
\begin{eqnarray}\label{e:pest2}
&&   \bP_y   \left( T\leq a_0 \delta (y)^2/4 \hbox{ and }
   Z^{\delta}_T   \in B(y_0 , \, \delta_0 r/16) \
 \Big|  \
   \sup_{s<T\wedge (a_0 \delta (y)^2/4)}
   |Z^\delta_{s}-y|\leq  \delta (y)/2  \right) \nonumber \\
 &\geq &
  \left( 1-e^{-\kappa_1 a_0 \delta(y)^2/4 } \right)
 \left(\kappa_2/\kappa_1 \right) \geq  c \, a_0 \, \delta (y)^2.
\end{eqnarray}
  By Meyer's construction \cite[Lemma 3.6]{BBCK} and Lemma \ref{L:exit2},
\begin{eqnarray*}
&& \bP_y \left( \sup_{s<T\wedge (a_0 \delta (y)^2/4)}
|Z^\delta_{s}-y|
\leq  \delta (y)/2 \right) \\
&\geq &  \bP_y \left(  \sup_{s\le a_0 \delta (y)^2/4}
|Z^\delta_{s}-y|\leq  \delta (y)/2 \hbox{ and } T\geq a_0 \delta (y)^2/4\right) \\
& \geq & e^{-\kappa \cdot a_0 \delta (y)^2/4}
\bP_y \left( \sup_{s\leq a_0  \delta (y)^2/4}
 |Z^{\delta}_s -y|\leq \delta (y)/2 \right)
 \geq 1/(4e^{\kappa}).
\end{eqnarray*}
This together with \eqref{e:pest2} yields that
\begin{equation}\label{e:pest3}
 \bP_y \left(T\leq a_0 \delta (y)^2/4,
 \  Z^{\delta}_T   \in B(y_0 , \, \delta_0 r/16) \hbox{ and } \sup_{s<T}
 |Z^\delta_{s_0}-y|\leq  \delta (y)/2  \right) \geq
  c\, \delta (y)^2.
\end{equation}
Since $s_0=t_1/3< a_0  (\delta_0 r)^2/16$, we have from
 Lemma \ref{L:exit2} that
$$ \inf_{x\in \R^d\setminus \NN} \bP_x \left( \sup_{   s\leq  s_0}
 |Z^\delta_s  -Z^\delta_0| \leq \delta_0 r/4 \,
  \right)
 \geq 1/4.
$$
Therefore we have by \eqref{e:pest} and \eqref{e:pest3} that
$$ \bP_y \left( Z^\delta_{s_0} \in B(0, (1-\delta_0/2) r) \hbox{ and } \tau_{B_r} > s_0 \right)
\geq c (r-|y|)^2 .
$$
Now for $t\in [t_1/2, t_2]$, $y\in B_r$ and $z\in B(0, (1-\delta_0/2) r)$,
by \eqref{e:ql}
\begin{eqnarray*}
q^{\delta, B_r}(t, y, z) &\geq &\int_{B(0, (1-\delta_0 /2) r)} q^{\delta,
B_r}(s_0, y, w)q^{\delta, B_r}(t-s_0, w, z) dw \\
&\geq & c \int_{B(0, (1-\delta_0/2) r)} q^{\delta, B_r}(s_0, y, w) dw \\
&=& c \, \bP_y \left( Z^\delta_{s_0} \in B(0, (1-\delta_0/2) r) \hbox{ and }
\tau_{B_r}    > s_0 \right)  \\
&\geq&  c (r-|y|)^2.
\end{eqnarray*}
This together with the Chapman-Kolmogorov's equation
$$ q^{\delta, B_r}(t, x, y) \geq \int_{B(0, (1-\delta_0/2) r)} q^{\delta, B_r}(t/2, x, z)
 q^{\delta, B_r}(t/2, z, y) dz
 $$
 proves the lemma.  \qed

\bigskip

\noindent {\bf Proof of Proposition \ref{thm:bbcklem47}.}~
(i) First, similarly to the proof of \cite[Lemma 4.1]{BBCK}, we have
\begin{equation}\label{41bbck}
q^{\delta,B_r}(t, x, y)\leq c_1
t^{-d/2}
\quad \hbox{ and } \quad
\left|\frac{\partial q^{\delta,B_r}(t, x, y)}{\partial t}\right|
\leq c_1 t^{-1-d/2}
\end{equation}
for every $x, y\in B_r$ and $t>0$. Using this,
$q^{\delta,B_r}(t,\cdot,y_0)\in \sF^{\delta,B_r}$ can be proved in the same
way as the proof of \cite[Lemma 4.5]{BBCK}.
Next, by Lemma \ref{thm:bbckpro43} and by the choice of $a_1$,
for every $y_0\in B_r$, $\varepsilon\in (0, 1)$ and
$\gamma \in \left(\frac{2-\beta}6, \, 1 \right]$,
there is a constant $C=C(y_0, \beta, \delta,\varepsilon)>0$ such that
\begin{equation}\label{44bbck}
\Psi_r(x)^\gamma/q^{\delta,B_r}(t,x, y_0)\leq C, \qquad
\hbox{for every }\, t\in (\varepsilon,\varepsilon^{-1}]
\hbox{ and } \,  x \in B_r. \end{equation}

Using this, $\Psi_r (\cdot )^{1/2}/q^{\delta,B_r}(t, \cdot, y_0)$ is
bounded on $B_r$. By extending the function
$x\mapsto \frac{\Psi_r (x)}{q^{\delta,B_r}(t, x, y_0)}$
to be zero on $B_r^c$, we see that
it vanishes continuously on $B_r^c$. Similar  to the proof of
Proposition \ref{P:charaF},
$$\sF^{\delta,B_r}= \left\{f\in L^2(\R^d; m_d):
f|_{B_r^c}\equiv 0 \hbox{ and }  \sE^\delta (f,f)<\infty \right\}.
$$
 So, in order to
prove $h_t(\cdot):=\Psi_r (\cdot )/q^{\delta,B_r}(t, \cdot, y_0)\in
\sF^{\delta,B_r}$, it is enough to prove $\sE^\delta
(h_t,h_t)<\infty$. Let $u_t(\cdot)= q^{\delta,B_r}(t, \cdot, y_0)$.
In order to show $\int_{B_r}\nabla h_t(x) A(x) \nabla h_t(x)
dx<\infty$, it is enough to prove
$\int_{B_r}|u_t (x)\nabla\Psi_r (x)-\Psi_r (x) \nabla u_t(x)|^2/ u_t (x)^4dx<\infty$,
since  $a(\cdot)$ is uniform elliptic. Computing this,
\begin{eqnarray*}
\int_{B_r}\frac{|u_t (x)\nabla\Psi_r (x)-\Psi_r (x)\nabla u_t (x)|^2}{u_t(x)^4}dx
 &\le & 2\left(\int_{B_r}\frac{|\nabla\Psi_r (x)|^2}{u_t (x)^2}dx+
\int_{B_r}\frac{|\Psi_r (x)\nabla u_t (x)|^2}{u_t (x)^4}dx\right)\\
&\le& 2 \left( c_1c_2^2 \, m_d(B_r)+c_2^4\int_{B_r}|\nabla u_t (x)|^2dx\right)<\infty,
\end{eqnarray*}
where $|\nabla\Psi_r|^2/\Psi_r\le c_1$ and $\Psi_r^{1/2}/u_t\le c_2$ (due to
(\ref{44bbck})) are used in the second
inequality. The proof of
$$ \int_{B_r}\int_{B_r}(u_t(x)-u_t(y))^2J_\delta(x,y)dxdy+ 2 \int_{B_r}u_t(x)^2
 \left(\int_{B_r^c}J_\delta(x,y)dy \right) dx<\infty
$$
 can be
done similarly to that of \cite[Lemma 4.6]{BBCK} (with a suitable change due
to the shape of $J_\delta$, for example $\gamma=(2-\beta)/3$ in the proof).
We thus obtain $\sE^\delta (h_t,h_t)<\infty$.

(ii) Given (i), (\ref{41bbck}) and (\ref{44bbck}), this
can be proved in the same
way as the proof of \cite[Lemma 4.7]{BBCK}. \qed

The idea of the proof of the following theorem is
motivated by that of Theorem 3.4 in \cite{CKK} and
Proposition 4.9 in \cite{BBCK}.
However, due to the existence of the divergence form part, various non-trivial
changes are required.

\begin{thm}\label{T:3.3}
 For each $t_0>0$, there exists $c =c (t_0)>0$, independent of
$\delta\in (0, 1)$ such that for every $x_0\in \R^d$,
$t\le t_0$,
 \begin{equation}\label{eqn:estball}
q^{\delta, B(x_0, t^{1/2})} (t, x,y)\ge c\, t^{-d/2}\qquad \hbox{ for q.e. }
x, y\in B(x_0, \sqrt{t}/2)
 \end{equation}
and
 \begin{equation}\label{eqn:est}
q^\delta (t, x,y)\ge c \, t^{-d/2}\qquad \hbox{ for q.e. }
x, y \hbox{ with } |x-y|^2 \leq t.
 \end{equation}
\end{thm}
\pf
Fix $\delta \in (0, 1)$ and, for simplicity, in this proof  we
sometimes drop the superscript ``$\delta$" from $Z^\delta$
and $q^\delta (t, x, y)$. Also, for notational convenience, let $x_0=0$.
For ball $B_r:=B(0,r)\subset \R^d$, let $q^{B_r}(t,x,y)$ denote the
transition density function of the subprocess $Z^{B_r}$ of $Z$ killed
on leaving the ball $B_r$.

Define $B:=B(0,1)$ and for $r\le 1$, let
$(\sE^{\angel{r}}, \sF^{\angel{r},B})$ be the Dirichlet form corresponding
to $\{r^{-1}Z^{\delta, B_r}_{r^2t}, t\geq
0\}$, which is the subprocess of $\{Z^{\angel{r}}_t:=r^{-1}Z^\delta_{r^2t},, t\geq
0\}$ killed on leaving the unit ball $B$. Define
\begin{equation}\label{eqn:3.8}
q_r^B(t, x, y)=q_r^{\delta,B}(t, x, y):=r^{d}q^{B_r}(r^2t, rx, ry).
\end{equation}
It is easy to see $q^B_r(t, x, y)$ is the transition density
function for process $r^{-1}Z^{\delta, B_r}_{r^2t}$.

Set $\Psi(x)=c((1-|x|)_+)^{a_1}$,
where $c>0$ is the normalizing constant.
Let $x_0\in B(0,1)$, $r\le 1$, and define
\begin{eqnarray*}
u(t, x)&:=& q_r^B(t, x, x_0), \\
v(t, x)&:=& q_r^B(t, x, x_0)/\Psi(x)^{1/2},\\
H(t)&:=& \int_{B} \Psi (y) \log u(t, y) dy , \\
G(t)&:=& \int_{B} \Psi (y) \log v(t, y) dy =\int_{B}
\Psi (y) \log u(t, y) dy
-\frac12 \int_{B} \Psi (x) \log \Psi (x) dx \\
&=& H(t) +c_1.
\end{eqnarray*}
By Proposition \ref{thm:bbcklem47} and the scaling, we have
\begin{equation}  \label{G-Eeq}
G'(t)  =- \sE^{\angel{r}} \Big( u(t, \cdot), \frac{\Psi}{u(t, \cdot)}
\Big)=:-(J_1+J_2),
\end{equation}
where $J_1$ is the diffusion part and $J_2$ is the jump part of the
Dirichlet form.

We first estimate the jump part. Write $J^{\angel{r}}_\delta(x, y):=r^{d+2} J_\delta (rx, r y)$.
By the same argument as in the proof of
Proposition 4.9 of \cite{BBCK}
(up to the formula fourth lines after (4.15) there), we have
\begin{eqnarray*}
&&J_2=\sE^{\angel{r},j} \Big( u(t, \cdot), \frac{\Psi}{u(t, \cdot)}\Big)\\
&\le & \int_B\int_B\{(\Psi(x)^{1/2}-\Psi(y)^{1/2})^2
-(\Psi(x)\wedge \Psi(y))(\log \frac{v(t,y)}{v(t,x)})^2\}J^{\angel{r}}_\delta(x,y)dxdy\\
&&~~+\int_B\Psi(x)\Big(2\int_{B^c}J^{\angel{r}}_\delta(x,y)dy\Big)dx\\
&\le & \int_B\int_B(\Psi(x)^{1/2}-\Psi(y)^{1/2})^2J^{\angel{r}}_\delta(x,y)dxdy
+\int_B\Psi(x)\Big(2\int_{B^c}J^{\angel{r}}_\delta(x,y)dy\Big)dx\\
&=&\sE^{\angel{r},j}(\Psi^{1/2},\Psi^{1/2}) \\
&\le & c_2r^{2-\beta }\sE(\Psi^{1/2},\Psi^{1/2}) \\
&\le&  c_2\sE(\Psi^{1/2},\Psi^{1/2}) <\infty,
 \end{eqnarray*}
where the last inequality is due to the shape of $J$ and the Lipschitz continuity
of $\Psi$ (note that $c_2\sE(\Psi,\Psi)$ is independent of $r$).

We next estimate the diffusion part.
\begin{eqnarray}
&&J_1=\sE^{\angel{r},c} \Big( u(t, \cdot), \frac{\Psi}{u(t, \cdot)}\Big)
\le \int_B\nabla u(t,x)a(rx)\nabla\Big(\frac{\Psi(x)}{u(t,x)}\Big)dx\nonumber\\
&= & \int_B\nabla\log u(t,x)a(rx)\nabla \Psi(x)dx-
\int_B\nabla\log u(t,x)a(rx)\nabla\log u(t,x)\Psi(x)dx.\label{wjidbvk}
\end{eqnarray}
Note that
\begin{eqnarray*}
0&\le &\int_B((\nabla \log u)\sqrt{\Psi}-\frac{\nabla\Psi}{\sqrt{\Psi}})
a^{\angel{r}}\cdot ((\nabla \log u)\sqrt{\Psi}-\frac{\nabla\Psi}{\sqrt{\Psi}})dx\\
&=& \int_B\nabla\Psi a^{\angel{r}}\cdot \nabla \Psi \Psi^{-1}dx+\int_B
(\nabla \log u) a^{\angel{r}}\cdot (\nabla \log u)\Psi dx-2\int_B(\nabla \log u) a^{\angel{r}}
\cdot\nabla \Psi dx,
\end{eqnarray*}
where $a^{\angel{r}}(\cdot)=a(r\cdot)$.
Using this and (\ref{unielli}) in (\ref{wjidbvk}), we obtain
\[J_1\le  c_3\int_B
\frac{|\nabla\Psi(x)|^2}{\Psi(x)}dx-c_4\int_B|\nabla\log u(t,x)|^2\Psi(x)dx=
c_5-c_4\int_B|\nabla\log u(t,x)|^2\Psi(x)dx,\]
where the last equality is due to the fact
$|\nabla\Psi(x)|^2/\Psi(x)\le c_{5.5}$ for $x\in B$, which is because
$a_1\ge 2$ in the definition of $\Psi$. Thus,
using Proposition \ref{wtedPI},
\[J_1\le c_6-c_7\int_B(\log u(t,x)-H(t))^2\Psi(x)dx.\]
Combining these, we obtain from (\ref{G-Eeq}),
\begin{equation}\label{firststep}
G'(t) =H'(t)\geq
-c_8+c_{7} \int_B (\log u(t,y)-H(t))^2 \Psi(y) \, dy.
\end{equation}
Given this inequality, (\ref{e:p1}) and Lemma \ref{L:exit2},
the rest of the proof is the same as that of
\cite[Proposition 4.9]{BBCK} (cf. also  \cite[Theorem 3.4]{CKK}).

\bigskip

\noindent {\bf Proof of Theorem \ref{T:3.4}.}
For any ball $B\subset \R^d$, let $ (\sE^{\delta, B}, \sF^{\delta,
B})$ denote the Dirichlet form of the subprocess $Z^{\delta, B}$ of
$Z^\delta$ killed upon leaving the ball $B$.
 Similarly to the proof of \cite[Theorem 1.5
and Theorem 2.6]{BBCK}, we can show that $(\sE^\delta, \sF^\delta)$ and
$(\sE^{\delta, B}, \sF^{\delta, B})$ converge as $\delta \to 0$ to
$(\sE, \sF)$ and $(\sE^B, \sF^B)$,
respectively  in the sense of
Mosco, where $B$ is a ball in $\R^d$. Therefore the semigroup of
$Z^\delta$ and $Z^{\delta, B}$ converge in $L^2$ to that of $X$ and
$X^B$, respectively. Theorem \ref{T:3.4} follows from Theorem \ref{T:3.3}
by a similar argument as that   for \cite[Theorem
1.3]{BBCK}. \qed

\section{H\"older continuity and Parabolic Harnack inequality}

\subsection{H\"older continuity}

 In this subsection, the Dirichlet form $(\sE, \sF)$ is given by \eqref{e:DF}
 with the jumping kernel satisfying the conditions \eqref{e:J1}-\eqref{e:J2}, and $X$
 is its associated strong Markov process in $\R^d$.

For  $r\in (0, 1]$, define
$$
  Q( x,r) :=(0,   r^2]\times B(x,r).
$$
For each $A\subset [0, \infty)\times \R^d $,
denote $\sigma_A:=\inf\{t> 0: Z_t\in A\}$
and $A_s:= \{y\in \R^d: \ (s, y)\in A\}$.

\begin{lem}\label{L:exit3}
There exists $C_2>0$ such that for all $x\in \R^d $,
 $0<r\leq 1$
and any compact subset $A\subset Q(0,x,r)$,
\[
 \bP^{(r^2, x)}
(\sigma_A<\tau_r)\ge C_2
\frac {m_{d+1} (A)}{r^{d+2}},\]
where $\tau_r=\tau_{Q(  x, r)}$ and
$m_{d+1}$ is the Lebesgue measure on $\R^{d+1}$.
\end{lem}

\pf
For $0<r\leq 1$,
\begin{eqnarray*}
  r^2 \, \bP^{(r^2, x)} (\sigma_A<\tau_r)
  &\ge & \int_0^{r^2} \bP^{(r^2,x)} \left((r^2-s,X_s^{B(x,r)})\in A\right)ds\\
&=& \int_0^{  r^2} \int_{A_{r^2-s}}p^{B(x,r)}(s,x,y)dyds\\
&\ge&
\int_0^{  r^2} \int_{A_{r^2-s}}\frac c{r^d}dyds=c \, \frac {m_{d+1} (A)}{r^d},
\end{eqnarray*}
where Theorem \ref{T:3.4} is used in the last inequality.
  \qed

\medskip

We can now establish the H\"older continuity for parabolic functions of $X$.
First, recall
the following  well-known formula (see, for example \cite[Appendix A]{CK2}).

\begin{lem}\label{L:4.9}(L\'evy system formula)
Let $f$ be a non-negative measurable function on ${\mathbb R}_+
\times \R^d\times \R^d$ that vanishes along the diagonal. Then for every
$t\geq 0 $, $x\in \R^d\setminus \NN$ and  stopping time $T$ (with
respect to the filtration of $X$),
$$
\E_x \left[\sum_{s\le T} f(s,X_{s-}, X_s) \right]=
\E_x \left[ \int_0^T \left( \int_{\R^d} f(s,X_s, y) J(X_s,y) dy
\right) ds \right].
 $$
\end{lem}

\medskip

\noindent{\bf Proof of Theorem \ref{T:holder}.}
For $x\in \R^d\setminus \NN$ and $r< 1$, apply
Lemma \ref{L:4.9} to
$f(s, y, z)= \1_{B(x, r)} (y) \1_{B(x, 2r)}(z)$
and $T=\tau_{B(x, r)}$.
Then it follows from \eqref{e:J2} and
Lemma \ref{L:exit},
for every $s\geq 2r$,
\begin{eqnarray}
 \bP_x \left( X_{\tau_{B(x, r)}} \notin B(x, s) \right)
 & = & \E_x \left[ \int_0^{\tau_{B(x, r)}} \left( \int_{\R^d \setminus B(x, s)}
   J(X_t, y) dy \right) dt \right] \nonumber \\
 &\leq &  4 (s\wedge 1)^{-2}\, \E_x \left[ \int_0^{\tau_{B(x, r)}} \left( \int_{\R^d }
  (|X_t-y|^2 \wedge 1 )  J(X_t, y) dy \right) dt \right] \nonumber \\
 &\leq & c \, (s\wedge 1)^{-2}\,  \E_x \left[ \tau_{B(x, r)}\right] \nonumber \\
 &\leq & c \, r^2/(s\wedge 1)^2.  \label{e:hitting}
\end{eqnarray}
Using this and Lemma \ref{L:exit3}, the rest of the proof is the same
as that for  the proof of Theorem 4.14
in \cite{CK1} except that the estimate for
\begin{equation}\label{eq:ck1forml}
\sum_{i=1}^\infty  \E_{z_1} \left[ q(Z_{\tau_{k+1} })-q(z_2); \,
           \sigma_A> \tau_{k+1} \hbox{ and } Z_{\tau_{k+1} }\in Q_{k-i}
      \setminus  Q_{k+1-i} \right]
\end{equation}
at the bottom of page 57 of \cite{CK1} should be bound as follows.
Take $\rho <\eta$, then
\begin{eqnarray*}
(\ref{eq:ck1forml})&\leq &
    \sum_{i=1}^k(b_{k-i}-a_{k-i}) \bP_{z_1}
   ( X_{\tau_{k+1} } \notin Q_{k+1-i}) + \| h \|_{\infty, R} \bP_{z_1}
   (X_{\tau_{k+1}} \notin Q_0)  \\
&\leq &   \sum_{i=1}^k  c  \, \eta^k  (\rho^2
    / \eta )^i  + c\, \| h \|_{\infty, R}
   \, \rho^{k+1} \\
&\leq & c\,  \eta^{k-1} \rho^2 + c\,  \rho^{k+1} \\
&\leq & c \eta^{k+1}.
\end{eqnarray*}
        \qed

\subsection{Parabolic Harnack inequality}
In this subsection, the Dirichlet form $(\sE, \sF)$ is given by \eqref{e:DF}
 with the jumping kernel satisfying the conditions \eqref{e:J1}-\eqref{e:J2} and \eqref{e:J3}, and $X$
is its associated strong Markov process in $\R^d$.

Recall that $Z_s:=(V_s, X_s)$ is the
 space-time process of $X$, where $V_s=V_0- s$.
The following lemma corresponds to \cite[Lemma 4.2]{CKK}.
Noting that the continuous component of the process does
not play any role since the function $h$ is supported
in $[0,\infty)\times B(x_0,3R)^c$, the proof is almost the same as that of
\cite[Lemma 4.2]{CKK}.
We point out that condition \eqref{e:J3} is used in a crucial way
in the proof of this lemma.

\begin{lem}\label{L:compa}
Let $R\leq 1$ and $\delta <1$.
$Q_1=[t_0+2\delta R^2/3, \, t_0+ 5\delta R^2]\times B(x_0,
3R/2)$, $Q_2=[t_0+ \delta R^2/3, \, t_0+11\delta R^2/2]\times
B(x_0, 2R)$ and define
$Q_-$ and $Q_+$ as in Theorem \ref{T:PHI}. Let $h: [0,\infty)\times
\R^d\to \R_+$ be bounded and supported in $[0,\infty)\times
B(x_0,3R)^c$. Then there exists $C_1=C_1(\delta)>0$ such that the
following holds:
\[\E^{(t_1,y_1)}[h(  Z_{\tau_{Q_1}})]\le C_1\E^{(t_2,y_2)}
[h(  Z_{\tau_{Q_2}})] \qquad \mbox{for }  (t_1,y_1)\in Q_- \mbox{ and
}  (t_2,y_2)\in Q_+.\]
\end{lem}

\medskip

\noindent{\bf Proof of Theorem \ref{T:PHI}.}
With the above lemma, Lemma \ref{L:exit3} and the heat kernel
estimates in the previous sections, the proof is almost the same as
that of the proof of \cite[Theorem 4.1]{CKK} for $R\le 1$. \qed

\section{Heat kernel upper bound estimate
under condition  \eqref{e:J4}}

For the remaining two sections, we assume that the jumping kernel $J$ for
the Dirichlet form $(\sE, \sF)$ of \eqref{e:DF} satisfies condition \eqref{e:J4}.
For simplicity,  define
$$\tp(r):=r^2\wedge \phi(r).
$$
Note that $r\to \tp (r)$ is a strictly increasing function on $[0,
\infty)$ so it has an inverse function $\tp^{-1} (r)$. Clearly,
$$   \tp^{-1} (r) = r^{1/2}\vee \phi^{-1}(r),
$$
where $\phi^{-1}$ is the inverse function of $\phi$.
 Note that
$$  \tp^{-1} (t)^{-d}= t^{-d/2}\wedge\phi^{-1}(t)^{-d} .
$$

\begin{thm}\label{upperHK}
There are positive constants $c_1$ and $ c_2$ such that for every
$x,y\in \R^d$ and $t>0$, we have
\begin{equation}\label{eqn:heat}
p(t, x, y) \leq c_1\, \tp^{-1} (t)^{-d} \wedge
\left( p^c(t, c_2|x-y|)+p^j(t, |x -y|) \right).
\end{equation}
\end{thm}

\bigskip

Before proving this theorem, we make some preparations.
For $r>0$,
let $(\sE^{\angel{r}}, \sF^{\angel{r}})$ be the Dirichlet form corresponding
to $\left\{ X^{\angel{r}}_t:=r^{-1}X_{\tp(r)t}, t\geq 0 \right\}$.
By simple computations, we see that
$\sF^{\angel{r}}=W^{1,2}(\R^d)$ and
for $u,v\in \sF^{\angel{r}}$,
\[\sE^{\angel{r}} (u,v)=\frac{\tp(r)}{r^2}\int_{\mathbb R^d}
\nabla u(x)\cdot a(rx)\nabla v(x)dx+\int_{\mathbb R^d}
(u(x)-u(y))(v(x)-v(y))J^{\angel{r}}(x,y)dxdy,\]
where $J^{\angel{r}}(x,y)=\tp(r)r^dJ(rx,ry)$. Note that
\[J^{\angel{r}}(x,y)\asymp \frac{\tp(r)}{|x-y|^d\phi(r|x-y|)}
=\frac 1{|x-y|^d\phi_r(|x-y|)},\]
where $\phi_r(s):=\phi(rs)/\tp(r)$ (note that $\phi_r$ enjoys
the properties (\ref{polycon}) and (\ref{polycon2})
with the constant $c>0$ independent of $r$).
Clearly the transition density function $p_r (t, x, y)$ of $X^{\angel{r}}$ with
respect to $m_d$ is given by
\begin{equation}\label{eqn:003.8}
p_r(t, x, y):=r^{d}p(\tp(r)t, rx, ry).
\end{equation}

The following on-diagonal estimate holds for $p(t, x, y)$:
\begin{equation}\label{eqn:2-0}
 p(t, x, y) \leq c\,
 \left( t^{-d/2}\wedge \phi^{-1}(t)^{-d} \right),\qquad\forall t>0.
 \end{equation}
If follows from the Nash inequality for the stable-type Dirichlet
form obtained in \cite[Theorem 3.1]{CK2}, we have $p(t, x, y) \leq
c \phi^{-1}(t)^{-d}$, so that (\ref{eqn:2-0}) holds. Thus, using
(\ref{eqn:003.8}), we have
\begin{equation}\label{eqn:2}
 p_r(t, x, y) \leq r^d  \left( \tp^{-1} (\tp (r) t) \right)^{-d}=:g(r, t).
 \end{equation}
 Clearly $g(r, 1)=1$ and
 \begin{eqnarray*}
 g(r, t)
 &\leq & c \left( r^d (\tp (r) t)^{-d/2} \1_{\{\tp(r)t\le 1\}}
  + r^d (\phi^{-1} ( \tp (r) t))^{-d}  \1_{\{\tp(r)t > 1\}} \right)\nonumber \\
 &\leq& c\,
 \left(  r^{d}\, \tp (r)^{-d/2} t^{-d/2} \1_{\{\tp(r)t\le 1 \}}
 + r^d \, \tp (r)^{-d/\beta_2} t^{-d/\beta_2} \1_{\{\tp(r)t> 1 \}} \right).
 \end{eqnarray*}

For $\lambda >0$, define
$$ J^{\angel{r, \lambda}} (x, y):= J^{\angel{r}}(x, y) \1_{\{|x-y|\leq \lambda \}}
$$
and let $(\sE^{\angel{r, \lambda}}, W^{1,2}(\R^d))$ be defined as $(\sE^{\angel{r}}, \sF^{\angel{r}})$
but with jumping kernel $J^{\angel{r, \lambda}}$ in place of $J^{\angel{r}}$.
Let $X^{\angel{r, \lambda}}$ be the symmetric strong Markov process
associated with $(\sE^{\angel{r, \lambda}}, W^{1,2}(\R^d))$.
The process $X^{\angel{r, \lambda}}$ can be obtained from $X^{\angel{r}}$ by removing
all the jumps whose size is larger than $\lambda$.
 We will apply  Davies' method to derive heat kernel
estimate for process $X^{\angel{r,\lambda}}$.
 On-diagonal estimate (\ref{eqn:2}) together with Theorem 3.25 of
\cite{CKS} implies that there exist constants $C>0$ and $c>0$,
independent of $\lambda >0$ and $\delta>0$ such that
\begin{equation}\label{eqn:kern4}
p_r^{(\lambda)}(t,x,y)\leq g(r,t) \;\exp\left(
-|\psi(y)-\psi(x)|+C\;
\Lambda_{r,\lam}(\psi)^2\; t\right)
\end{equation}
for all $t>0$,  $x,y\in \R^d\setminus \NN$ and every
$\lambda>0$, and for some $\psi$ satisfying $\Lambda_{r,\lam}(\psi)<\infty$,
where
\begin{equation*}
\Lambda_{r,\lam}(\psi)^2= \|e^{-2\psi}\Gamma_{r,\lambda}[e^\psi]\|_\infty \vee
\|e^{2\psi}\Gamma_{r,\lambda}[e^{-\psi}]\|_\infty.
\end{equation*}
Here
\begin{equation}  \label{eqn:densi}
\Gamma_{r,\lambda}[v](\xi)=\frac{\tp(r)}{r^2}
\sum_{i,j=1}^da_{ij}(r\xi)\frac{\partial v}{\partial x_i}
(\xi)\frac{\partial v}{\partial x_j}(\xi)+
\int_{|\eta-\xi|\leq \lambda}
(v(\eta)-v(\xi))^2J^{\angel{r}}(\eta,\xi) d \eta, ~~\xi\in
\R^d.
\end{equation}
Define
$$
\mathcal{H}(\Gamma_{r,\lambda}):=\left\{v: G\to \mathbb{R}\; \Big|
\;\sup_{\xi\in \R^d}\Gamma_{r,\lambda}[v](\xi)<\infty\right\}.
$$
A key observation is that $\mathcal{H}(\Gamma_{r,\lambda})$ contains
the \textit{cut-off distance function} $\psi$ given by
\begin{equation}\label{eqn:testf}
\psi(\xi):=\frac s3
\left(  |\xi-x| \wedge |x- y| \right)
\qquad \hbox{for } \xi \in \R^d,
\end{equation}
where $s>0$ is a parameter to be chosen later. Note that
$|\psi(\eta)-\psi(\xi)|\leq (s/3)
 |\eta-\xi|$ for all $\xi,\eta\in \R^d$. So
\begin{align*}
e^{-2\psi(\xi)}\Gamma_{r,\lambda}[e^\psi](\xi)&\le
c_1|\nabla \psi (\xi)|^2+
\int_{|\eta-\xi|\leq \lambda}
(1-e^{\psi(\eta)-\psi(\xi)})^2J^{\angel{r}}(\eta,\xi)
d \eta\\
&\leq c_1\frac{s^2}9+\int_{|\eta-\xi|\leq \lambda}
(\psi(\eta)-\psi(\xi))^2\;
e^{2|\psi(\eta)-\psi(\xi)|}J^{\angel{r}}(\eta,\xi) d \eta\\
&\leq c_1\frac{s^2}9+
(\frac s3)^2\, e^{2 s\lambda/3} \int_{|\eta-\xi|\leq
\lambda}|\eta-\xi|^2
J^{\angel{r}}(\eta,\xi) d\eta\\
&\leq c_1\frac{s^2}9+c s^2e^{2 s\lambda/3}\int_0^\lambda
\frac {t}{\phi_r (t)} \, dt \\
&\leq c_1\frac{s^2}9+c s^2e^{2 s\lambda/3}
\frac{\lambda^2}{\phi_r (\lambda)}\\
&\leq c_2(s^2+
\frac{e^{s\lambda}}{\phi_r (\lambda)}),
\end{align*}
for every $\xi\in \R^d$. Here we used Lemma 2.1(ii) of \cite{CK2}
 for the fourth inequality and
the fifth inequality is by (\ref{polycon2}).
The same estimate holds for
$e^{2\psi(\xi)}\Gamma_{r,\lambda} [e^{-\psi}](\xi)$.
Denote the  constant $c_2>0$ by $C_*$ and define
\begin{equation}\label{eqn:5-0}
F(  r, \lambda, s, t, R):=\exp \left(-\frac {sR}3+ C_*\left( s^2+\frac {e^{s\lambda}}
{\phi_r (\lambda)} \right)t
\right).
\end{equation}
Then, by (\ref{eqn:kern4}), with $R=|x-y|$, we have
\begin{equation}\label{eqn:chan1230}
p_r^{(\lambda)}(t, x, y) \leq
g(r,t) F(  r, \lambda, s, t, R).
\end{equation}
Note that there is a  freedom to choose $s>0$ properly. We are now
ready to prove Theorem \ref{upperHK}.

\medskip

\noindent {\bf Proof of Theorem \ref{upperHK}.} \
By (\ref{eqn:2-0}), it suffices to show that
\begin{equation}\label{eqn:3}
p(t, x, y) \leq c_1\,\left(  p^c(t, c_2 |x- y|)+p^j(t, |x- y|)
\right).
 \end{equation}
Our proof consists of considering 5 cases. Recall that $R:=|x-y|$.

{\bf Case 1:} $R^2<t<\phi (R) \leq 1$.

Take $r=1$, $\lambda=R$ and
$s=\frac{1}{\sqrt{t}}$ in \eqref{eqn:chan1230}. Note that in this case,
$g(1, t) = ct^{-d/2}$ and
$$ \frac{e^{sR}}{\phi (R)} = \frac{ e^{R/\sqrt{t}}}{\phi (R)}
< \frac{e}{t} =es^2.
$$
So
 \begin{eqnarray*}
  p_1^{(\lambda)}(t, x, y)
   \leq c_1\,  t^{-d/2}   e^{-\frac{sR}3 +C_*(1+e)s^2 t}
=  c_2\,  t^{-d/2} \,  e^{-\frac R{3\sqrt{t}} } .
\end{eqnarray*}
(In fact, $p_1^{(\lambda)}(t, x, y) \leq c_1 t^{-d/2}$ in this case.)
It follows by Meyer's construction that
\begin{eqnarray*}
 p(t, x, y) &\leq&  p_1^{(\lambda)} (t, x, y) + t\, \sup_{x\in \R^d}
\int_{\R^d} J(x , y)   \1_{\{ |x -y|>\lambda\}} dy \\
 &\leq& c_1\,  t^{-d/2} \,  e^{-R/(3\sqrt{t}) } + c_1
\frac{t}{R^d \phi (R)} \\
&\leq &  c_1\,  t^{-d/2} \, e^{-R^2/(3t) } + c_1
\frac{t}{R^d \phi (R)}.
\end{eqnarray*}
The last inequality is due to the assumption that $R^2< t$.
So \eqref{eqn:3} holds in this case.

\bigskip

{\bf Case 2:} $\phi (R) \leq t$.

This is a free lunch as $ p^j (t, x, y) \approx c \phi^{-1}
(t)^{-d}$ in this case and \eqref{eqn:3} follows.

\bigskip

Let $K=\beta_1/(72C_*(d+\beta_1))$ and let $a=eK/c$,  where $C_*$ and $c$ are  the
positive constants in
\eqref{eqn:5-0} and (\ref{polycon}), respectively.
Before we consider the remaining three cases, let us first do estimate
on $F:=F(r, \lambda, s, t, R)$ under two situations:
$$ \hbox{ (i) } \    e^{KR^2/t}\geq \frac{a\phi_r (R)}{t} \hbox{ with }R^2\geq t,
\qquad   \hbox{ and } \qquad \hbox{ (ii) } \
e^{KR^2/t}< \frac{a\phi_r (R)}{t}.
$$
Since $\displaystyle \min_{x>0}   e^x/x= e$, we have
\[\frac 1K\cdot\frac t{\phi_r(R)}e^{KR^2/t}=
\frac {R^2}{\phi_r(R)}\cdot\frac t{KR^2}e^{KR^2/t}\ge \frac{\tp(r)R^2}{\phi(rR)}
\cdot e ,\]
which, by \eqref{polycon} is no less than $1/c$ if $\min\{r, R\}\geq 1$ or
if $r\leq 1$ but $rR\geq 1$.
So Situation  (ii) may happen only when $r<1\leq R$ and $rR<1$.

 \medskip

\noindent  {\bf Situation (i):}
  $e^{KR^2/t}\geq \frac{a\phi_r (R)}{t}$ and $R^2\geq t$.

Let $H=\beta_1/(12(d+\beta_1))$. We   take   $\lambda=HR$ and $s=(HR)^{-1}
\log (e\phi_r (R)/t)>0$ in \eqref{eqn:5-0}.
By \eqref{polycon}, there is a constant $c_1>0$
such that
\[\frac {e^{s\lambda}}{\phi_r (\lambda)}t
\le c_1 \, \frac {e^{s\lambda}}{\phi_r (R)}\, t =c_1e.
\]
Moreover, using the assumption,
\begin{eqnarray*}
C_*s^2t&=&C_*\frac {st}{HR}\log \frac{e\phi_r(R)}t
=C_*\frac {st}{HR}\log \frac ea +C_*\frac {st}{HR}
\log \frac{a\phi_r(R)}t\\
&\le& c_2\frac {st}R + C_* \frac{st}{HR} \frac {KR^2}{t}
=s \left( c_2\frac tR+\frac {R}6 \right)\le \frac {sR}4+ c_3,
\end{eqnarray*}
since $K=\beta_1/(72C_*(d+\beta_1))=H/(6C_*)$. The last inequality is
due to that fact that when $R^2/t\ge 12c_2$,
$$ s \left( c_2\frac tR+\frac {R}6 \right) \leq s \left( \frac{R}{12} +
  \frac {R}6 \right) =  \frac {sR}4,
  $$
 while for  $1\le R^2/t< 12c_2$,
\[c_2\frac{st}R=c_2\frac t{HR^2}\log \left( e\frac{\phi_r (R)}t \right)
\le \frac{c_2}H\log \left(\frac ea e^{12c_2K} \right)=:c_3.
\]
So, by (\ref{eqn:5-0}), we have
\begin{equation}\label{e:s1}
F\le \exp \left( -\frac {sR}{12}+c_3+C_*c_1e \right) =c_4
\left(\frac t{\phi_r (R)e} \right)^{1/(12H)} =c_5 \left(\frac
t{\phi_r (R)}\right)^{d/\beta_1+1}.
\end{equation}

\noindent  {\bf Situation (ii):}
 $e^{KR^2/t}< \frac{a\phi_r (R)}{t}$.

We take  $\lambda=KR/(6C_*)$,
$s=R/(6C_* t)$ in \eqref{eqn:5-0}.
By \eqref{polycon}, there is a constant $c>0$
such that
$$
\frac{e^{s\lambda}}{\phi_r (\lambda)}t
\le  c\, \frac{e^{s\lambda}}{\phi_r (R)} \, t =
 c\, \frac{e^{KR^2/t}}{\phi_r (R)} \, t \le ca .
$$
  So
\begin{eqnarray}
F &\le & \exp \left(-\frac{sR}3+C_*s^2t+C_*ca \right) \label{e:s2} \\
&=& c_6 \exp \left(-\frac{sR}3+C_*\frac {sR}{6C_*} \right)
=c_6\exp \left(-\frac{sR}6 \right) =c_6 \exp \left( -\frac{R^2}{6C_*t} \right). \nonumber
\end{eqnarray}

{\bf Case 3:} $t\leq 1\leq R$.

We will take $r=1$ in this case so by \eqref{eqn:2},
$$ g(t, 1)= c t^{-d/2} \leq c t^{-d/\beta_1}.
$$
This case falls into Situation (i) and so we have from
\eqref{eqn:chan1230} and \eqref{e:s1}
$$ p_1^{(\lambda)}(t, x, y) \leq ct^{-d/\beta_1}
\left(\frac t{\phi (R)} \right)^{d/\beta_1+1}=
c\frac t{\phi (R)^{d/\beta_1+1}}\leq
  \frac {c_7t}{R^d\phi (R)},$$
where we used (\ref{polycon})
in the last inequality. By Meyer's construction,
we conclude
\begin{equation}\label{scalehkca3}
 p (t, x, y) \leq
 c_8  \, \Big(\frac{t}{R^d \, \phi_r (R)}+
 \frac t{R^d\phi(R)}\Big)
 \le   \, \frac{c_8 t}{R^d \, \phi (R)}.
\end{equation}
This establishes \eqref{eqn:3} in this case.

\bigskip

{\bf Case 4:} $\phi (R) \geq t\geq 1$.

Let $r=\phi^{-1}(t)\ge 1$, $x'=x/r$ and $y'=y/r$. Since $R\ge r$,
$|x'-y'|\ge 1$ so the estimate for $p_r (1, x', y')$
 falls into Situation (i). As $g(r, 1)=1$, we have from
\eqref{eqn:chan1230}, \eqref{e:s1}  and Meyer's construction
\begin{eqnarray*}
r^{d}p(\phi(r), x, y)&=& p_r (1, x', y') \\
&\leq & p^{(\lambda)}_r (1, x', y') + \sup_{x\in \R^d} \int_{\R^d}
 J^{\angel{r}} (x, y) \1_{\{|x-y|>\lambda\}} \\
&\leq &   c
\left(\frac 1{\phi_r (|x'-y'|)} \right)^{d/\beta_1+1} + \frac{c}{|x'-y'|^d \phi_r (|x'-y'|)} \\
&\leq &  c_9  \, \frac{1}{\phi_r(1)^{d/\beta_1} |x'-y'|^d \, \phi_r (|x'-y'|)}
    + \frac{c_9}{|x'-y'|^d \phi_r (|x'-y'|)} \\
& \le& \frac{c_{10}\phi (r)}{|x'-y'|^d \phi_r (|x-y|)}.
 \end{eqnarray*}
 Here we used \eqref{polycon} in the second to the last inequality and the fact
 that $\phi_r (1) \geq 1$ in the last inequality.
Since $t=\phi(r)$, we conclude that
\[p(t, x, y)\le \frac{c_{10} t}{|x-y|^d \, \phi(|x-y|)}.\]
This proves \eqref{eqn:3} in this case.

\bigskip

{\bf Case 5:} $t<R^2(\le \phi (R)) \le 1$.

Let $r=R=|x-y|$, $x'=x/r$, $y'=y/r$. Note that $\tp (r)=r^2$ as $r\leq 1$ and
$|x'-y'|= 1$. Let $t'=t/r^2\le 1$.  Note that
$$ g(r, t') \leq c (t')^{-d/2} \leq c (t')^{-d/\beta_1}.
$$

If $e^{K /t'}\ge  a\phi_r (1)/t'$, then we are in Situation (i)
for $p_r (t', x', y')$.
By the same calculation as that for Case 3, we have
\[r^{d}p(r^2t', x, y)=p_r (t', x', y')
 \le  \frac{c_{11} t'}{|x'-y'|^d \, \phi_r (|x'-y'|)}
 =\frac{c_{11}t'r^2}{|x'-y'|^d \, \phi(|x-y|)}.\]
Noting $t=t'r^2$, we obtain
\[p(t, x, y)\le \frac{c_{11} t}{|x-y|^d \, \phi(|x-y|)}.\]

If $e^{K/t'}< a\phi_r (1)/t'$,
then we are in Situation (ii) for
$p_r(t', x', y')$. So by \eqref{eqn:chan1230}, \eqref{e:s2} and
Meyer's construction
\begin{eqnarray*}
r^{d}p(r^2t', x, y) &=&p_r (t', x', y') \\
&\leq & p_r (t', x', y') + t' \sup_{x\in \R^d} \int_{\R^d}
 J^{\angel{r}} (x, y) \1_{\{|x-y|>\lambda\}} dy \\
&\leq& c_{12} t'^{-d/2} \exp \left(-\frac{c_{13} |x'-y'|^2}{t'}\right)
 + \frac{c_{14} t'}{|x'-y'|^d\phi_r(|x'-y'|)}.
\end{eqnarray*}
Noting $t=t'r^2$, we obtain
\[p(t, x, y)\le c_{15} t^{-d/2} \exp (-\frac{c_{16} |x-y|^2}t) +
\frac{c_{17}t}{|x-y|^d \, \phi(|x-y|)}.\]
This proves the claim (\ref{eqn:3}).

The upper bound estimate in
 (\ref{eqn:heat}) is now established for every $t>0$ and $x, y\in \R^d$. \qed

\section{Heat kernel lower bound estimate under condition \eqref{e:J4}}

Recall that  $\tp (t):=t^2 \wedge \phi (t)$ and so
$\tp^{-1} (t)^{-d} =t^{-d/2}\wedge\phi^{-1}(t)^{-d}$.
In this section, we will establish the following.

\begin{thm}\label{lowerHK}
There exist positive constants $c_1$ and $c_2$ such that

\begin{equation}\label{eqn:heatlb}
p(t, x, y) \geq c_1\,
 \tp^{-1} (t)^{-d}
\wedge \left( p^c(t, c_2 |x-y|)+p^j(t, |x- y|) \right)
\end{equation}
for each $x,y\in \R^d$ and $t>0$.
\end{thm}

To prove it, we need first establish some tightness results
and extend Lemma \ref{L:exit3} to all $r>0$ and Theorem \ref{T:PHI}
to all $R>0$.

\subsection{Tightness and some lower bound estimate}

Using the heat kernel upper bound, we can prove the following
estimate of the exit time from a ball.
\begin{propn}\label{P:4.2}
For each $A> 0$ and $0<B<1$, there exists $\gamma =\gamma (A,B)\in
(0,\, 1/2)$ such that for every $r>0$
and $x\in \R^d\setminus \NN$,
\[\bP_x \left( \tau_{B(x,\, Ar)}<\gamma \, \tp(r) \right)
\le B.\]
\end{propn}
\pf Let $x\in \R^d\setminus \NN$. By the upper bound estimate in
(\ref{eqn:heat}), for every $s>0$ and $t>0$,
\begin{eqnarray*}
\bP_x \left(|X_t-x| \geq s\right)
  &=& \int_{B(x,s)^c} p
  (t, x, y) dy\\
 &\le & \int_{B(x,s)^c} \frac{c_1\, t dy}{|x- y|^d\phi (c_1|x-y|)}
 +c_2t^{-d/2}\int_{B(x,s)^c}\exp (-\frac{c_3|x-y|^2}t) dy\\
 &\le &\frac {c_4t}{\phi (s)}+c_5\exp (-\frac{c_6s^2}t)
 \le \frac {c_4t}{\phi (s)}+\frac {c_7t}{s^2}\le\frac {c_8t}{\tp (s)}.
\end{eqnarray*}
The above computation is standard; see Lemma 2.1(i) in \cite{CK2}
for the estimate of the stable part in the second inequality, and
\cite{B98} Lemma 3.9 (a)
for the estimate of the Gaussian part in the second inequality.
Given this inequality, the rest of the proof is the same as that of
Proposition 4.9 in \cite{CK2}
 with $\tp$ in place of $\phi$ for the case of $\gamma_1=\gamma_2=0$ there.
\qed

\medskip

Using Proposition \ref{P:4.2}, one can prove the following proposition
in the same way as the proof of Proposition 4.11 in \cite{CK2}
but with $\tp$ in place of $\phi$ for the case of $\gamma_1=\gamma_2=0$ there.

\begin{propn}\label{lower2}
There exist constants
 $c_1 \ge 2$ and $c_2>0$ such that for every $t>0$ and
 every $x, y\in \R^d\setminus \NN$ with
\begin{equation}\label{hbound19}
\bP_{x} \left(X_t\in B(y, c_1\tp^{-1}(t)) \right)
  \geq c_2\frac{t(\tp^{-1}(t))^d}{|x-y|^d\tp ( |x-y|)}.
    \end{equation}
\end{propn}

\subsection{Parabolic Harnack Inequality}

Denote $\gamma (1/2, 1/2)$ in Proposition \ref{P:4.2} by $\gamma_0$.
For each $r,t>0$, we define
\[   Q(t,x,r) :=[t,t+\gamma_0 \tp(r)]\times B(x,r).
\]
The following is an extension of Lemma \ref{L:compa} to all $r>0$.

\begin{lem}\label{L:4.10}
There exists $C_1 >0$ such that for every
 $x\in \R^d$, $r>0$, $y\in B(x, \, r/3)$ and
a bounded nonnegative function $h$ on $[0, \infty)\times \R^d $
that is supported in $[0,\infty)\times B(x,2r)^c$,
\begin{equation}\label{compah}
\E^{(\gamma_0 \tp (r), x)}
\left[ h(\tau_r,X_{\tau_r}) \right]\le C_1
  \E^{(\gamma_0 \tp (r), y)}
\left[ h(\tau_r,X_{\tau_r}) \right],
\end{equation}
where $\tau_r=\tau_{Q(0,x,r)}$.
\end{lem}

\pf The proof is the same as Lemma 6.1 in \cite{CK2}.
Note that the continuous component of the process does
not play any role since the function $h$ is supported
in $[0,\infty)\times B(x,2r)^c$.
(Note that in \cite{CK2} the space-time process is running forward
in the sense that $V_t=V_0+t$ there while in this paper $V_t=V_0-t$
is defined to run backward. Clearly there is one-to-one correspondence
between these two situations. Thus the estimate
in Lemma 6.1 in \cite{CK2} is
under probability law $\bP^{(0, x)}$ while here it is under
$\bP^{(\gamma_0 \tp (r), x)}$.
The same remark applies in the following when \cite{CK2} is cited, for example,
in the proof of the next three results.)
\qed

\medskip

For each $A\subset [0, \infty)\times \R^d $,
denote $\sigma_A:=\inf\{t> 0:   Z_t\in A\}$.

\begin{lem}\label{L:4.12}
There exists $C_2>0$ such that for all $x\in \R^d $, $r>0$ and
any compact subset $A\subset Q(0,x,r)$,
\[
 \bP^{(\gamma_0 \tp (r), x)}
(\sigma_A<\tau_r)\ge C_2
\frac {m_{d+1} (A)}{r^d\tp(r)},\]
where $\tau_r=\tau_{Q(0, x, r)}$.
\end{lem}

\pf
When $r\leq 1$, this is proved in Lemma \ref{L:exit3}.
 When $r\ge 1$, we have $\tp(r)=\phi(r)$ so the desired inequality
can be proved similarly to Lemma 6.2 in \cite{CK2}.
\qed

\medskip

Define
$U(t,x,r):=\{t\}\times B(x,r)$.

\begin{cor}\label{C:4.13}
For every $0<\delta \le \gamma_0$,
there exists $C_3>0$ such that for
every $R\in (0,1]$, $r\in (0,  R/4]$ and $(t,x)\in Q(0,z,R/3)$ with
$0<t\leq \gamma_0 \tp (R/3)- \delta \tp (r)$,
$$
 \bP^{(\gamma_0 \tp (R/3), z)}
\left(\sigma_{U(t,x,r)}<\tau_{Q(0,z,R)} \right)\ge
   C_3 \frac{r^d\tp(r)}{R^d\tp(R)} .
$$\end{cor}
\pf Given Lemma \ref{L:4.12} and Proposition \ref{P:4.2},
the proof is the same as Corollary 6.3 in \cite{CK2}
but with $\tp$ in place of $\phi$ there.
\qed

\medskip

The following extends the parabolic Harnack principle in Theorem \ref{T:PHI}
to all $R>0$.

\begin{thm}\label{P:4.4}
For every $0<\delta \le \gamma_0$,
there exists $c_1>0$
such that for every $z\in \R^d$, $R>0$ and
every non-negative function $h$ on $[0, \infty)\times \R^d$
that is parabolic and bounded on
$[0,\gamma \tp(2R)]\times B(z,2R)$,
\[\sup_{(t,y)\in Q(\delta \tp(R),z,R)}h(t,y)\le c_1 \,
\inf_{y\in B(z,R)}h(0,y).\]
In particular, the following holds for
$t>0$.
\begin{equation}\label{harpt}
\sup_{(s,y)\in
Q((1-\gamma)t,z,\tp^{-1}(t))}
p(s,x,y)\le c \,
\inf_{y\in B(z, \tp^{-1}(t))}p((1+\gamma)t,x,y).
\end{equation}
\end{thm}
\pf Given Lemma \ref{L:4.10}, Lemma \ref{L:4.12} and
Corollary \ref{C:4.13}, the proof of this PHI is
the same as that of Theorem 4.12 in \cite{CK2}
(see also the proof of Theorem 4.5 in \cite{SV07}).
\qed
\subsection{Lower bound}
\begin{lem}\label{L:4.7}
There exist $c_1, \, c_2>0$ such that
\[p(t, x, y)\ge c_1 \, (\tp^{-1}(t))^{-d}\]
for all $t>0$ and
$x,y\in \R^d\setminus \NN$ with $|x-y|\le c_2 \, \tp^{-1}(t)$.
\end{lem}
\pf This is already proved in Theorem \ref{T:3.4} for $t\le 1$.
Given (\ref{eqn:heat}), Proposition \ref{P:4.2}, and
Theorem \ref{P:4.4}, the proof is the same as that of
Lemma 4.13 in \cite{CK2}
but with $\tp$ in place of $\phi$ there.
\qed

\bigskip

\noindent
{\bf Proof of Theorem \ref{lowerHK}}.
Let $t>0$.
Due to Lemma \ref{L:4.7}, it is enough to prove
the theorem for $|x-y|\ge c_2 \, \tp^{-1}(t)$.
Applying Proposition \ref{lower2} with $t_* =(1-\gamma) t$
in place of $t$, we have
\[\bP_x  (X_{t_* }\in B(y, \, c_1\tp^{-1}(t_*)))\ge c_2
\frac{t_*(\tp^{-1}(t_*))^d}{|x-y|^d\phi (c_3
|x-y|)}.
\]
 As $m_d (B(y, \, c_1 \phi^{-1}(t_*)))\le c_4 (\phi^{-1}(t_*))^d$, the above
implies $p (t_*, x, z)\ge c_5 \, t/(|x-y|^d\phi(c_3|x-y|))$
for some $z\in B(y, \, c_1 \phi^{-1}(t_*))$. By applying
(\ref{harpt}) as before, we have
\[p(t,x,y)\ge c\frac t{|x-y|^d\phi(|x-y|)}.\]
For (\ref{eqn:heatlb}), the exponential decay
appears on the RHS only when $t<r^2(\le \phi (r)) \le 1$ (Case 4 in the upper bound), where
$r=|x-y|$.
So, the only case left is this case. In this
case, choose $N\in \bN$ so that
$s:=t/N\asymp (r/N)^2$ (so $N\asymp r^2/t$).
Then, $p(s,x,y)\ge c s^{-d/2}$,
by Lemma \ref{L:4.7}. Thus the usual chain argument
gives $p(t,x,y)\ge ct^{-d/2}\exp (-c'r^2/t)$. \qed

\vskip 0.3truein

\noindent {\bf Zhen-Qing Chen}

\noindent Department of Mathematics, University of Washington,
Seattle, WA 98195, USA

\noindent E-mail: zchen@math.washington.edu

\bigskip

\noindent {\bf Takashi Kumagai}

\noindent Department of Mathematics, Faculty of Science, Kyoto
University, Kyoto 606-8502, Japan.

\noindent E-mail: kumagai@math.kyoto-u.ac.jp

\end{document}